\documentclass[9pt,shortpaper,twoside,web]{IEEEtran}
\usepackage{cite}

\usepackage{url}

\ifCLASSINFOpdf
\usepackage[pdftex]{graphicx}
\usepackage{epstopdf}
\usepackage{subfigure}
\graphicspath{{./}{./figures/}}
\else
\fi
\usepackage{xcolor}

\usepackage{amsmath}
\usepackage{amssymb}

\newtheorem{property}{Property}
\newtheorem{nnassumption}{\bf Assumption}

\newtheorem{nntheorem}{\bf Theorem}
\newenvironment{thm}{\begin{nntheorem}\it}{\end{nntheorem}}
\newtheorem{nncorollary}{\bf Corollary}

\newtheorem{nndefinition}{\bf Definition}

\newtheorem{nnproposition}{\bf Proposition}
\newenvironment{prop}{\begin{nnproposition}\it}{\end{nnproposition}}
\newtheorem{nnproblem}{\bf Problem}

\newtheorem{nnlemma}{\bf Lemma}
\newenvironment{lemma}{\begin{nnlemma}\it}{\end{nnlemma}}
\newtheorem{nnremark}{\bf Remark}
\newenvironment{rem}{\begin{nnremark} \rm }{\hfill \hspace*{1pt}\hfill $\circ$\end{nnremark}}
\newtheorem{nnexample}{\bf Example}

\newenvironment{proof}{{\bf Proof.}}{\hfill \hspace*{1pt}\hfill $\Box$}

\allowdisplaybreaks

\begin{document}
		%
		\title{Observer-based output feedback for an age-structured SIRD model}
		%
		%
		%
		
		\author{Candy Sonveaux, Christophe Prieur, Gildas Besançon, Joseph J. Winkin
			\thanks{Candy Sonveaux is with University of Namur, Department of Mathematics and naXys, Rue de Bruxelles 61, B-5000 Namur, Belgium (e-mail: candy.sonveaux@unamur.be).\newline\indent Gildas Besançon and Christophe Prieur are with Univ. Grenoble Alpes, CNRS, Grenoble INP, GIPSA-lab, 38000 Grenoble, France.\newline\indent Joseph J. Winkin is with 
				University of Namur, Department of Mathematics and naXys, Rue de Bruxelles 61, B-5000 Namur, Belgium.\newline\indent 				
				The first author wishes to thank the FRS-FNRS as well as the PhC Tournesol Program that supported research stays to GIPSA-Lab, supporting the collaboration between UNamur (Belgium) and Gipsa-Lab (France).
			}
		}
		
		%
		%

	\markboth{}%
	{Sonveaux \MakeLowercase{\textit{et al.}}}
	%



	\maketitle
	
	\begin{abstract}
		 An age-structured Susceptible-Infected-Recovered-Deceased (SIRD) epidemic model is considered. The aim of this paper is to design an observer-based output feedback control law, representing an immunization process, typically  vaccination, intended to decrease the peak of infected individuals in the population. At first, well-posedness and stability of the system in open-loop are investigated. Then, to obtain the observer-based output feedback law, a state feedback law is designed by using a normal form. Conditions to ensure stability are established. However, due to physical constraints, this law needs to be adapted. Therefore, a constrained state-feedback law is implemented. This law is designed to fulfill the physical constraints while having good properties (Lipschitz for instance), needed for the last part of the article. Finally, an observer-based output feedback law is obtained using high-gain observer. At each step of the design, convergence properties are obtained. Finally, numerical simulations are performed.    
	\end{abstract}
	

	%
	\IEEEpeerreviewmaketitle
	
		\section{Introduction}
\indent In order to face an epidemic, it is of major importance to understand the underlying dynamical processes. Once the dynamics of an epidemic is understood, we may want to act on it. This can be done by treating infected people or by enhancing immunity before \textcolor{black}{individuals are infected by} the disease (typically by vaccination). In the second approach, an interesting question that arises is ``how to implement an immunization strategy by taking into account the ages of the individuals?". This question is particularly important because several factors in disease propagation are age-dependent, see \cite{dowd_2020} for instance. This is the case for childhood diseases such as measles (e.g. \cite{Tudor_1985}), the AIDS epidemic (e.g. \cite{Griffiths_2000}) and also for cattle diseases (e.g. \cite{Brock_2020}), where cattle management is influenced by the age of the animals. \\\indent This paper aims to address this previous essential question. This question has been tackled for instance in \cite{Griffiths_2000}, in the context of the AIDS epidemic, drawing conclusions based on HIV age profiles and common sense. In the present work, a systematic strategy based on control theory is developed. This methodology has been used by various researchers (see \cite{Hernandez_2020} or \cite{Ram} in the case of COVID-19 pandemic) but not always applied to age-structured models with output feedback. \\ \indent To address this question, a model considering the age of individuals is needed. The model that is used is an extension of the well-known SIR (Susceptible-Infected-Recovered) model of \cite{kermack_1991}, where an additional category is added, the group $D$ of dead individuals. Moreover, the age is a key factor of this research. Therefore the SIR model is adapted into an age-structured SIRD model. There are two ways to describe age-dependent epidemic models, either by considering the age to be continuous (see for instance \cite{Inaba1990}, \cite{Hethcote_2008}, \cite{Jia} or \cite{sonveaux}), or by discretizing it (as done for instance in \cite{andreasen_frommelt_2005}). In this paper, the second option is used. Therefore each individual lays in a given class of age and in a given state of health (susceptible, infected, recovered or dead). \textcolor{black}{Hence, the age represents the age of an individual, contrary to the work in \cite{Webb} for instance, where the infection age is considered.} The dynamics of the disease propagation is described by a set of $4n$ nonlinear ordinary differential equations (ODE), where $n$ denotes the number of classes of age. The model considered here can be seen as an age-discretized version of the partial integro-differential equation model studied in \cite{sonveaux}, where the aging effect is not considered, since it is not relevant for short-term diseases. \\ \indent In order to solve the question of interest, an  immunization strategy, modeled under the form of a feedback control approach can be adopted.  Some studies on the impact of some control policy for fatal illness, using ODE models, were done (e.g. as in \cite{calafiore_2022}). But given that the precise full state of the system is typically unknown in practice, output feedback strategies become necessary. This issue has been explored concerning epidemic processes in \cite{Niazi_2023}, along with an optimal control method. Here an alternative observer-based output feedback control is presented with the purpose to decrease the peak of infected people. \\ \indent First, the designed control law is a linearizing feedback for the model in normal-form that reduces the maximum number of infected individuals. Conditions ensuring the non-negativity of the input are established. Besides the fact that the latter is mandatory due to its physical interpretation, it is shown that, for the case of one class of age, such an input yields a lower peak of infected individuals with respect to the zero-input case. \\ \indent Then, due to feasibility constraints (limited number of actions per unit of time), this work presents an amplitude limited feedback for which asymptotic convergence is investigated.  However, this law cannot be designed in practice because it requires the knowledge of the whole state. Therefore a high-gain observer is designed in order to build an observer-based age-dependent output feedback law. The analysis of the closed-loop system resulting from observer-based feedback is based on a separation principle. It is shown that such a design leads to performance recovery, namely pointwise attractivity.\\ \indent
	The paper is organized as follows. In Section \ref{Sec_Model}, the age-structured SIRD model studied in this paper is presented. The dynamical analysis of the open-loop model is performed in Section \ref{Sec_dyn_an}. Results on boundedness of the states and on pointwise asymptotic stability are obtained. Then, two state feedback laws are designed in Section \ref{Sec_pos_state_feed}. The observer issue is then addressed in Section \ref{Sect_output_fdb} , leading to an observer-based output feedback law. Finally, numerical simulations illustrating the analytical results are performed in Section \ref{Sec_num_sim} on an example inspired by literature results on COVID-19 and real data.
	
	\section{Model Formulation}\label{Sec_Model}
	In this section, the age-dependent SIRD epidemic model used in the following is described. It is inspired by the papers \cite{calafiore_2022} and \cite{franco_2021}, but here an input representing an immunization process is added. \\ \indent In SIRD models, the population is divided in four distinct classes: the class $S$ of susceptible individuals who can get the disease, the class $I$ of infected individuals who can transmit the disease; those infected individuals can either recover from the disease and join the group $R$ of recovered individuals, or die and go to the $D$ class of dead individuals. The age of the individuals is taken into account by considering classes of ages. In the following, $n$ classes are considered. Therefore, the dynamics of the disease propagation is given by a set of $4n$ nonlinear ordinary differential equations 	\begin{equation}\label{SIRD}
		\left\{\begin{array}{rl}
			\dfrac{dS_k}{dt}\left(t\right) &= -\lambda S_k(t) \displaystyle\sum_{j=1}^n \dfrac{M_{kj}I_j(t)}{N_j} -p_k\theta_k\left(t\right)S_k\left(t\right)\vspace{0.2cm}\\ 
			\dfrac{dI_k}{dt}\left(t\right) &= \lambda S_k(t) \displaystyle\sum_{j=1}^n \dfrac{M_{kj}I_j(t)}{N_j} -\left(\gamma_{R_k}+\gamma_{D_k}\right)I_k\left(t\right)\vspace{0.2cm}\\ 
			\dfrac{dR_k}{dt}\left(t\right) &=\gamma_{R_k}I_k\left(t\right)+p_k\theta_k\left(t\right)S_k\left(t\right)\vspace{0.2cm}\\ 
			\dfrac{dD_k}{dt}\left(t\right) &=\gamma_{D_k}I_k\left(t\right),
		\end{array}
		\right.
	\end{equation}
	
	\noindent where, for $k=1,...,n$, the state variables $S_k$, $I_k$, $R_k$ and $D_k$ denote the number of susceptible, infected, recovered and deceased individuals in the $kth$ class of age, respectively. Remark that, when considering only one class of age, i.e $n=1$, we recover the classical SIRD model with no class of age. The sum terms correspond to the fact that a susceptible individual in the $k$th class of age can be infected when a contact occurs with an individual of any class of age. Moreover, in those equations, the aging effect is not considered since migration of individuals from a class of age $k$ to the $(k+1)$th class of age is not taken into account. Therefore, this model is best suited for short-term disease. In view of the equations, a closed population is considered. This means that the total number of individuals in the population (deceased individuals included) is constant. Indeed, by summing the equations for a given $k$, it is obvious that $\frac{dN_k}{dt}$ is identically 0 where $N_k=S_k+I_k+R_k+D_k$ is the total number of individuals in the $kth$ class of age. Hence $N_k$ is constant. In this case, only $3n$ equations are needed because, by knowing $N_k=S_k(0)+I_k(0)+R_k(0)+D_k(0)$ and $3n$ appropriately chosen state variables, the last $n$ ones can be recovered. Moreover, the model parameters $M_{kj}\textcolor{black}{\in\left[0,N_j\right]}$, $\gamma_{R_k}\textcolor{black}{\in[0,1]}$ and $\gamma_{D_k}\textcolor{black}{\in[0,1]}$ are age-dependent and represent the social contact (average contact per unit of time from an individual of the $kth$ class of age with individuals of the $jth$ class), the recovery rate and the death rate, respectively. The parameter $\lambda$ denotes the transmission probability of the disease. Those parameters need to be calibrated thanks to parameters identification approaches. Concerning the coefficients $M_{kj}$, these can be derived through the use of the online tool Socrates 
		(Social contact rates) from \cite{social_contact_data}, \cite{willem_2020} where statistical surveys have gathered data about the contact habits of the participants. Finally, the function $\theta_k\left(t\right)$ is the input variable representing the rate of action on susceptible individuals in the class of age $k$ at time $t$. The parameter $p_k\in\left(0,1\right]$ represents the probability for the immunization to work for the $kth$ class of age. This immunization can, for instance, be seen as an ideal vaccination. Indeed, it is used on the susceptible individuals. So it is a preventive mechanism. However, model \eqref{SIRD} \textcolor{black}{relies on the assumption that the immunization process works} perfectly: once immunized, an individual never catches the disease. Improvements in this area require further study. 
	\section{Dynamical Analysis}\label{Sec_dyn_an}
	This section is devoted to the analysis of some properties of the open-loop model. More precisely, the fact that the state variables are non-negative and bounded is first underlined, showing consistency with the physical meaning of the problem. Indeed, the states should stay in the set \begin{align*}
			B&:=\left\{(S_1...S_n, I_1...I_n,R_1...R_n,D_1... D_n )^T \in \mathbb{R}^{4n}:\right.\\
		&\left. 0\leq S_k \leq N_k, 0\leq I_k \leq N_k, 0\leq R_k \leq N_k,   0\leq D_k \leq N_k \right. 
		\\  &\left.\text{for } k=1,...,n\right\}
	\end{align*} Moreover, a proof of the pointwise attractivity of the equilibria of the model in open-loop is also provided.\\ \\
	\textcolor{black}{In the next proposition, the space $L_{\text{loc}}^{\infty}(\mathbb{R}_+)$ is introduced. It describes the space of locally essentially bounded functions on $\mathbb{R}_+$, i.e $L_{\text{loc}}^{\infty}(\mathbb{R}_+)=\left\{f:\mathbb{R}_+\to \mathbb{R} \text{ where }\forall \text{ compact subset } K \subset \mathbb{R}_+, \right.\\ \left. \text{there exists a constant }C_K \text{ such that } \vert f(x)\vert\leq C_K \text{ almost }\right.\\ \left.\text{everywhere in } K\right\}$.\\ } 
	\begin{prop}\label{prop_pos_state}
	Let $\theta_k\in L_{\text{loc}}^{\infty}(\mathbb{R}_+)$ such that for $k=1,...,n$ and for all $t\geq 0$, $\theta_k(t)\geq 0$. Then the compact set $B$ is positively invariant for model \eqref{SIRD}.
	\end{prop}
	\vspace{0.3cm}
The proof of this proposition is based on the concept of essential non-negativity. Additional details can be found in \cite{sonveaux_thesis}.
	This proposition implies that the model is well-defined for dealing with numbers of individuals. Moreover, notice that this proposition applies to model \eqref{SIRD} in closed-loop with non-negative feedback law in $L_{loc}^{\infty}(\mathbb{R}^+)$. This fact is used in Theorem \ref{thm_pos_unsaturated_fdb}. Moreover, notice that this regularity condition on $\theta_k$, for $k=1,...,n$ ensures the existence of a unique global solution to the model \eqref{SIRD}, for all $t\geq 0$ and $x_0\in B.$ This is a consequence of \cite[Theorem 3.3]{khalil_2002}.\\ Another result concerning pointwise attractivity of the equilibria in open-loop can be obtained. A set $X$ is called pointwise attractive if every solution  of the dynamical system converges to one equilibrium of the set of equilibria, i.e. $\displaystyle\lim_{t\to\infty}x(t)$ exists in $X$. This is inspired by the concept of pointwise stability (see for instance \cite[Definition 3.1]{Hui_2008}),
but where the equilibrium points are not necessarily Lyapunov stable. Indeed, this requirement is not relevant in the study of epidemic models.\\
	\begin{prop}\label{prop_asympt_stab}
	For $k=1,...,n$, let $\theta_k\in L_{\text{loc}}^{\infty}(\mathbb{R}_+)$ be a control of feedback-type, $\theta_k(t)=\Theta_k(x(t))
		$ where $\Theta_k:B\to \mathbb{R}$ is a locally Lipschitz function 
		such that for all $x\in B$, $\Theta_k(x)\geq 0$. The set of disease-free equilibria of model \eqref{SIRD} $$X^{\star}=\left\{x\in B:I_k=0 \textcolor{black}{\text{ and } \theta_k S_k=0}, k=1,...,n\right\}$$ is pointwise attractive. In particular, any trajectory starting in $B$ converges to a given equilibrium point in this set, i.e. for all initial conditions $x_0=$ \small{$(S_1(0)...S_n(0), I_1(0)...I_n(0), R_1(0)...R_n(0), D_1(0)...D_n(0))^T\in $}
		\normalsize $B,$ there exists a unique equilibrium $  x^{\star}\in X^{\star}$ such that $x(t)\to x^{\star}$, as $t$ tends to infinity, where $x(t)$ is the solution of model \eqref{SIRD} with initial condition $x_0$.
	\end{prop}\vspace{0.3cm}
	\begin{proof}
		In what follows, $x^{\star}\in X^{\star}$ denotes any equilibrium of model \eqref{SIRD} and the corresponding input is denoted by $\theta^{\star}=\left(\theta_1^{\star},...,\theta_n^{\star}\right)$. Therefore, it follows from the last ODE of \eqref{SIRD} that $I_k^{\star}$ equals $0$ for $k=1,...,n$.  
		Moreover, the first and third equations imply to solve the relation $\theta_k^{\star}S_k^{\star}=0.$ Three cases can be considered, either $S_k^{\star}\textcolor{black}{\in (0,N_k]}$ then $\theta_k^{\star}=0$ or $S_k^{\star}=0$ and $\theta_k^{\star}\neq 0 $ or $S_k^{\star}=\theta_k^{\star}=0$. 
		Furthermore, since there is no condition on the other variables \textcolor{black}{and since they may take any value in $\left[0,N_k\right]$}, it follows that there is an infinity of possible equilibrium points $(S_1^{\star},...,S_n^{\star}, I_1^{\star}=0,...,I_n^{\star}=0, R_1^{\star},...,R_n^{\star}, D_1^{\star},..., D_n^{\star})$ corresponding to the disease-free case. 
		\\
		LaSalle's theorem (see for instance \cite[Theorem 4.4]{khalil_2002}) leads to the conclusion about the attractivity of the set $X^{\star}$. Indeed, considering $v=(S_1,...,S_n,I_1,...,I_n,R_1,...,R_n)$ such that $(v,D_1,...,D_n)\in B\subset \mathbb{R}^{4n}$ where $B$ is a compact set which is positively invariant with respect to model \eqref{SIRD} by Proposition \ref{prop_pos_state}, we can define the function $V(v)=\displaystyle\sum_{k=1}^{n} \left(S_k+I_k+R_k\right)\geq 0.$ The time derivative of $V$ along the trajectories of model \eqref{SIRD} is given by
		\begin{align*}
			\dot{V}(v)&=\displaystyle\sum_{k=1}^n\dfrac{dS_k}{dt}+\dfrac{dI_k}{dt}+\dfrac{dR_k}{dt}\\
			&=\displaystyle\sum_{k=1}^n -\gamma_{D_k}I_k\left(t\right)\\
			&\leq 0
		\end{align*}
		on $B$. Therefore, by LaSalle's theorem, any solution starting in $B$ converges to the set $X^{\star}$ as time goes to infinity.  Moreover, since $\dfrac{dS_k(t)}{dt}\leq 0$  and $S_k(t)\geq 0$ for all $t\geq 0$, it follows that $S_k(t)$ converges to an equilibrium point $S_k^{\star}$. The same applies to $R_k(t)$ that is increasing and upper bounded by $N_k$ for all $t\geq 0$. Therefore it converges to an equilibrium point $R_k^{\star}$. It follows that $D_k=N_k-S_k-R_k-I_k$ converges to $D_k^{\star}=N_k-S_k^{\star}-R_k^{\star}$ Finally, we can conclude the proof by combining LaSalle's theorem and convergence of bounded monotone functions which implies pointwise attractivity of the set. 
	\end{proof}
	\section{Non-negative bounded state feedback design}\label{Sec_pos_state_feed}
	The aim of this section is to design a control law to decrease the maximal total number of infected individuals in the population over time. This will reduce the impact of the disease in terms of hospital burden. The next proposition emphasizes that a control action can help to decrease the peak of infected individuals in the population. Then, two problems are solved in the following. The first one aims to develop a control law to get the asymptotic convergence of the state trajectories towards a disease free equilibrium with an exponential rate of convergence of the infected individuals towards zero. The other problem concerns the design of an amplitude limited control where the amplitude constraint is set a priori, in view of the maximum number of people who can get handled per unit of time. \\
	\begin{prop}\label{prop_goal_vacc}
		Consider two trajectories $x(t)$ and $\bar{x}(t)$ of model \eqref{SIRD} with k=1 \textcolor{black}{(hence the indices are not written)}, with the same initial conditions $x(0)=\bar{x}(0) \in B$ such that  \textcolor{black}{\begin{align}
			I(0)&=\bar{I}(0)\geq 1,\label{I0}\\
			S(0)&=\bar{S}(0)>\dfrac{N\left(\gamma_{R}+\gamma_{D}\right)}{\lambda M}:=\dfrac{\gamma}{\alpha}.\label{S0}
		\end{align}} The first trajectory $x(\cdot)$ is obtained without control, i.e. such that $\theta(t)=0$ for all $t\geq 0$, whereas the second one, $\bar{x}(\cdot)$, is such that  \textcolor{black}{$\exists$ $\delta>0$ : $$ \theta(t)> 0 \hspace{0.2cm}\forall t\in\left[0,\delta\right] \text{ and } \theta(t)\geq 0 \hspace{0.2cm}\forall t>\delta.$$}
		Then, the maximum number of infected individuals, all ages considered subject to the input $\theta(t)$ is smaller than the one without input, i.e $$\displaystyle\max_{t\geq 0} \bar{I}(t)< \displaystyle\max_{t\geq 0} I(t).$$
	\end{prop} \vspace{0.3cm}
	\begin{proof}
When only one class of age is considered, the dynamics of the infected individuals are given by $\dot{I}(t)=\left(\alpha S(t)-\gamma\right) I(t)$ and the dynamics for the susceptible are given by $\dot{S}(t)=\left(-\alpha I(t)-p_1\theta(t)\right)S(t)$.\\ \textcolor{black}{Thanks to \eqref{I0} and \eqref{S0}, $\dot{I}(0)>0$. Therefore, the trajectories for the infected individuals are initially increasing. Furthermore, the functions $I(t)$ and $\bar{I}(t)$ tend to $0$ as time tends to infinity, by Proposition \ref{prop_asympt_stab}. Hence, there exists time $T$ and $\bar{T}$ where $I(T)$ and $\bar{I}(\bar{T})$ equals $I(0)$ and $\bar{I}(0)$ respectively. Since the functions $I(t)$ and $\bar{I}(t)$ are continuous on $\left[0,T\right]$ and $\left[0,\bar{T}\right]$, it follows that $\displaystyle\max_{t\in \left[0,T\right]} I(t)$ and $\displaystyle\max_{t\in\left[0,\bar{T}\right]} \bar{I}(t)$} are reached at times $t_m$ and $\bar{t}_m$ such that $\dot{I}(t_m)=0$ and $\dot{\bar{I}}(\bar{t}_m)=0$ and $0< I(t_m), \bar{I}(\bar{t}_m)\leq N_1$. This implies that
$S(t_m)=\bar{S}(\bar{t}_m):=\hat{S}=\dfrac{\gamma}{\alpha}.$ \textcolor{black}{Moreover, after $t_m$ and $\bar{t}_m$, the susceptible trajectories remain lower than $\dfrac{\gamma}{\alpha}$ since they are decreasing. It follows that $\dot{I}(t)\leq 0$ and $\dot{\bar{I}}(\bar{t})\leq 0$ for all $t>t_m$ and $t>\bar{t}_m$ respectively. This implies that $\displaystyle\max_{t\in\left[0,T\right]} I(t)=\displaystyle\max_{t\geq 0} I(t)=:I_{\max}$ and $\displaystyle\max_{t\in\left[0,\bar{T}\right]} \bar{I}(t)=\displaystyle\max_{t\geq 0} \bar{I}(t)=:\bar{I}_{\max}$}\\ Now, considering $\dfrac{dS}{dI}=\dfrac{-\alpha S I - p_1 \theta S}{\alpha SI-\gamma I}$ \textcolor{black}{(where the set $\left\{(S,I):\alpha SI-\gamma I=0\right\}$ is of Lebesgue measure zero)} and integrating, one gets \begin{align*}
		\displaystyle\int_{S_0}^{\hat{S}} \dfrac{\alpha S-\gamma}{S} dS &=-\displaystyle\int_{I_0}^{\bar{I}_{\max}} \dfrac{\alpha I + p_1 \tilde{\theta}(I)}{I} dI\\
		\alpha\left(\hat{S}-S_0\right)-\gamma\ln\left(\dfrac{\hat{S}}{S_0}\right)&=-\alpha\left(\bar{I}_{\max}-I_0\right)-p_1\displaystyle\int_{I_0}^{\bar{I}_{\max}}\dfrac{\tilde{\theta}(I)}{I}dI,
		\end{align*}where $x(0)=\bar{x}(0)$ is used \textcolor{black}{and $\tilde{\theta}(I):=\theta(t(I))$, for $t\in \left[0,t_m\right]$. This can be defined since $I$ is an injective function (strictly increasing, $\dot{I}(t)>0$) for all $t<t_m$.} It follows that 
		\begin{equation}\label{eq_SI1}
		\bar{I}_{\max}=-\dfrac{\gamma}{\alpha}+S_0 +\dfrac{\gamma}{\alpha} \ln\left(\dfrac{\gamma}{\alpha S_0}\right)+I_0-\dfrac{p_1}{\alpha}\displaystyle\int_{I_0}^{\bar{I}_{\max}}\dfrac{\tilde{\theta}(I)}{I}dI,
		\end{equation}  In the case where no input is considered, \eqref{eq_SI1} still holds with $\tilde{\theta}$ replaced by $0$.  
		Hence, \begin{align*}
		I_{max}&=\bar{I}_{\max}+\dfrac{p_1}{\alpha}\displaystyle\int_{I_0}^{\bar{I}_{\max}}\dfrac{\tilde{\theta}(I)}{I}dI\label{eq_I_no_input}\geq \bar{I}_{\max}
		\end{align*} because $\bar{I}_{max}\geq I_0$, $\tilde{\theta}(I)> 0$ \textcolor{black}{for the first} $I$ and $0<I_0<I<\bar{I}_{\max}.$ 
	\end{proof}\\ \\
	Remark \textcolor{black}{that conditions \eqref{I0} and \eqref{S0} are trivially satisfied since there are automatically met at the beginning of an epidemic.} Moreover, this proposition implies that in order to decrease the maximum number of individuals in the whole population seen as a single class of age, it is better to act than doing nothing.\\ The extension of Proposition \ref{prop_goal_vacc} to the case of an arbitrary number of classes of age is still an open question.
	
	\subsection{Unconstrained Feedback law}
	In this subsection, the model is put into normal form to implement a control law that implies that the infected individuals for each class of age are bounded by a decreasing exponential function, improving the result obtained in Proposition \ref{prop_asympt_stab}, where only asymptotic stability is proved, without additional information about the rate of convergence.\\
	\subsubsection{Feedback design}\phantom{i}\\ \\
	In view of model \eqref{SIRD} and the control objectives, only the first $2n$ equations are needed in the feedback design. Indeed, for $k=1,...,n$, the control of $I_k(t)$ (and therefore $D_k(t)$) \textcolor{black}{relies only on} the dynamics of $S_k$ and $I_k$ \textcolor{black}{since the $(S_k,I_k)$-system is in cascade with the $(R_k,D_k)$-system}. Therefore, the $2n$ first equations of model \eqref{SIRD} are a nonlinear control affine system, written in state-space form as 
	\begin{subequations}
		\label{state_space}
		\begin{align}
			\hspace{-0.3cm}\left\{\begin{array}{lp{0.2cm}l}
				\dot{x}\left(t\right)&=&f\left(x\left(t\right)\right)+g\left(x\left(t\right)\right)\theta\left(t\right)\\
				y\left(t\right)&=&h\left(x\left(t\right)\right)
			\end{array}\right.
			\intertext{where $x\left(t\right)=\left[I_1\left(t\right),...,I_n\left(t\right), S_1\left(t\right),...,S_n\left(t\right)\right]^T\in \mathbb{R}^{2n},$ for all $ t \geq 0$ is the state space vector, $h\left(x\left(t\right)\right)=\left[I_1\left(t\right),...,I_n\left(t\right)\right]^T\in\mathbb{R}^{n}, $ for any $ t\geq 0$ is the output to be controlled, which is chosen equal to the infectious population and $\theta\left(t\right)=\left[\theta_1\left(t\right),...,\theta_n\left(t\right)\right]^T\in\mathbb{R}^{n},$ for any $t\geq 0$ is the input function. Moreover,}
			g\left(x\left(t\right)\right)=\begin{pmatrix}
				0_{n\times n}\\
				-p_k \text{diag}(S_k(t))_{k=1,...,n}
			\end{pmatrix}
			\intertext{where $\text{diag}(a_k)_{k=1,...,n}$ denotes the diagonal matrix whose elements are $a_1,...,a_n$, }\textit{} f\left(x\left(t\right)\right):=\left(
			f_1\left(x\left(t\right)\right) 
			...
			f_{2n}\left(x\left(t\right)\right)
			\right)^T
		\end{align}
	\end{subequations}
	where
	\begin{align}
		f_k\left(x\left(t\right)\right)&=
		\lambda	S_k\left(t\right)\displaystyle\sum_{j=1}^n C_{kj} I_j\left(t\right)-\left(\gamma_{R_k}+\gamma_{D_k}\right)I_k\left(t\right),\label{f_k}\\ 
		f_{n+k}\left(x\left(t\right)\right)&=-\lambda S_k\left(t\right)\displaystyle\sum_{j=1}^n C_{kj} I_j\left(t\right)\label{f_nk}
	\end{align} for $k=1,...,n$, where the notation $C_{kj}$ refers to $\dfrac{M_{kj}}{N_j}$ for readability.\\ 
	In this case, the relative degree of the system equals the dimension of the state space for any $x\in\mathcal{D}\subset \mathbb{R}^{2n}$ where \begin{align} \mathcal{D}=\Bigg\{x \in \mathbb{R}^{2n} \text{ s.t } S_k\neq 0\text{ and }\displaystyle\sum_{j=1}^n C_{kj} I_j\neq 0, k=1,...,n \Bigg\}.\label{D}
	\end{align} Hence, we can introduce a nonlinear invertible (for $x \in \mathcal{D}$) coordinate change 
	\begin{align}
		\label{change_of_variable}
		z_{1k}\left(t\right)&=h_k\left(x\left(t\right)\right)=I_k\left(t\right),\nonumber\\
		z_{2k}\left(t\right)&=L_fh_k\left(x\left(t\right)\right)=f_k\left(x\left(t\right)\right)
	\end{align}		for $k=1,...,n$. Notice that this change of variable is only invertible for $x \in \mathcal{D}$. However, the aim of the control is to bring the number of infected individuals to zero, which does not belong to $\mathcal{D}$. Hence, in Section \ref{Sub_saturated}, a constrained feedback law is implemented to take this fact into account.\\ Thanks to this change of variable, the model dynamics in normal form, in the neighborhood of any $x\in\mathcal{D}$ is given by 
	\begin{align}\label{normal_form}
		\hspace{-0.2cm}	\left\{\begin{array}{lp{0.2cm}l}
			\dfrac{dz_{1k}\left(t\right)}{dt}&=&z_{2k}\left(t\right),\\
			\dfrac{dz_{2k}\left(t\right)}{dt}&=&L_f^2h_k\left(x\left(t\right)\right)+L_{g_k}L_fh_k\left(x\left(t\right)\right)\theta_k\left(t\right)
		\end{array}
		\right.
	\end{align}
	for $k=1,...,n$, where \begin{align*}
		L_f^2h_k\left(x\left(t\right)\right)&=\lambda 
		S_k\left(t\right)\displaystyle\sum_{j=1}^n C_{kj} f_j\left(x\left(t\right)\right)\\ -\left(\gamma_{R_k}\right.&\left.+\gamma_{D_k}\right)f_k\left(x\left(t\right)\right)+\lambda f_{n+k}\left(x\left(t\right)\right)\displaystyle\sum_{j=1}^n C_{kj}I_j(t).
	\end{align*}
	Let $ 
	A\left(x\left(t\right)\right)=\text{diag}\left(-\lambda p_kS_k\left(t\right)\displaystyle\sum_{j=1}^n C_{kj}I_j\left(t\right)\right),$ $k=1,...,n,$
	$	v\left(x\left(t\right)\right)=\begin{pmatrix}
		v_1\left(x\left(t\right)\right) &\cdots & v_n\left(x\left(t\right)\right)
	\end{pmatrix}^T\label{v}$
	such that $ v_k\left(x\left(t\right)\right)=-\alpha_2^kf_k\left(x\left(t\right)\right)-\alpha_1^k I_k\left(t\right),$ where $\alpha_1^k$ and $\alpha_2^k$ are some free parameters that will be tuned in order to get stability with exponential convergence for the infected individuals. Moreover, $
	b\left(x\left(t\right)\right)=\begin{pmatrix}
		b_1\left(x\left(t\right)\right) & \cdots & b_n\left(x\left(t\right)\right)
	\end{pmatrix}^T,$
	where $ b_k\left(x\left(t\right)\right)=L_f^2h_k\left(x\left(t\right)\right)$. This allows us to design a stabilising and linearising state feedback, inspired by \cite[Chapter 5]{Isidori}. \\
	The state feedback control law defined by
		\begin{align}
			\theta\left(t\right)&=A^{-1}\left(x\left(t\right)\right)\left(v\left(x\left(t\right)\right)-b\left(x\left(t\right)\right)\right),\nonumber\\
				&=\left(
			\dfrac{v_1\left(x\left(t\right)\right)-L_f^2h_1\left(x\left(t\right)\right)}{L_{g_1}L_fh_1\left(x\left(t\right)\right)},
			\cdots,
			\dfrac{v_n\left(x\left(t\right)\right)-L_f^2h_n\left(x\left(t\right)\right)}{L_{g_n}L_fh_n\left(x\left(t\right)\right)}
		\right)^T\label{control}
			\end{align} for all $x(t)\in\mathcal{D}$, where $A$, $b$ and $v$ are defined above, applied on system \eqref{state_space}, induces the linear closed-loop dynamics given 
%
%
for $k=1,...,n$, by
		\begin{align}
			\left\{\begin{array}{lp{0.2cm}l}
				\dfrac{dz_{1k}\left(t\right)}{dt}&=&z_{2k}\left(t\right)\\
				\dfrac{dz_{2k}\left(t\right)}{dt}
				&=&-\alpha_2^kz_{2k}\left(t\right)-\alpha_1^kz_{1k}\left(t\right)	
			\end{array}\right.\label{CLdynamics}
		\end{align}
		which can be written as
		\begin{align*}
			\dot{z}\left(t\right)&=\begin{pmatrix}
				0_{n \times n} & Id_n\\
				-\tilde{A}_1 & -\tilde{A}_2
			\end{pmatrix}z\left(t\right),\\
			&:=\bar{A}z\left(t\right)
		\end{align*}
		where $z\left(t\right):=\left[z_{1_1}\left(t\right) \cdots z_{1_n}\left(t\right) z_{2_1}\left(t\right) \cdots z_{2_n}\left(t\right) \right]^T$.\\ \\
%
	Moreover, the infected individuals converge exponentially for the closed-loop system under some appropriate conditions.\\
	\begin{thm}\label{thm_stab_exp}
		Pick an initial condition $x_0 \in \mathcal{D}$, where $\mathcal{D}$ is given by \eqref{D}. Assume that the control tuning parameters satisfy $\alpha_j^k>0$ for $j=1,2$ and $k=1,...,n$. \\
		Then, as long as the closed-loop states  with the state feedback \eqref{control} remain in $\mathcal{D}$, the infected populations $I_k(t)$, for $k=1,...,n,$ of model \eqref{SIRD}, are bounded by a decreasing exponential function, i.e  there exist $
		\mu_k >0$ and $C_k >0$   such that $\forall$ $x_0\in \mathcal{D}$, $x(t)\in\mathcal{D}$ and $$I_k(t)\leq C_kI_k(0) e^{-\mu_k t}, \forall t\geq 0 \text{ and } k=1,...,n.$$ Moreover the susceptible, recovered and deceased populations converge asymptotically to some constants $S_k^{\star}, R_k^{\star}$ and $D_k^{\star}$ respectively,  for $k=1,...,n$.
	\end{thm} \vspace{0.3cm}
	\begin{proof}
		Since the closed-loop dynamics \eqref{CLdynamics} is a system of decoupled ODE's it can be written as
		\begin{align}
			\label{decoupled_ODE}
			\left\{
			\begin{array}{rll}
				\dot{z}_{new}\left(t\right)&=&\hat{A}z_{new}\left(t\right),\\
				y\left(t\right)&=&Hz_{new}\left(t\right)
			\end{array}
			\right.
		\end{align}
		with $z_{new}:=\begin{pmatrix}
			z_{1_1}&z_{2_1}&\cdots& z_{1_n}&z_{2_n}
		\end{pmatrix}^T$, $\hat{A}=\text{blockdiag}(\bar{A}_k)$, where $\bar{A}_k=\begin{pmatrix}
			0 & 1\\
			-\alpha_1^k & -\alpha_2^k
		\end{pmatrix}$ and $H=\begin{pmatrix}
			1 & 0 & 1 & \cdots & 1 & 0
		\end{pmatrix}.$ \\
		Therefore, $\hat{A}$ is stable if all its eigenvalues are in the open left half-plane. However, the eigenvalues of $\hat{A}$ are those of the $\bar{A}_k$'s matrices. 
		Since  $\alpha_1^k$ and $\alpha_2^k$ are positive, the real parts of the eigenvalues of $\bar{A}_k$ (hence of $\hat{A}$) are negative. Then the control law \eqref{control} exponentially stabilizes the model in normal form \eqref{normal_form}.	\\
		Therefore, $z_{new}\left(t\right)$ exponentially converges asymptotically to zero. It follows that $z_{1k}\left(t\right)=I_k\left(t\right)$ converges to zero as time goes to infinity for $k=1,...,n$. Moreover, by the convergence of bounded monotone functions, it follows that,  for $k=1,...,n$, $S_k(t)\to S_k^{\star}$, $R_k(t)\to R_k^{\star}$ and $D_k(t)\to D_k^{\star}$ as times goes to infinity. 
	\end{proof}\\ \\
Notice that the bound on $I_k(t)$ allows to find the following bound on $D_k^{\star}$,  $$ D_k^{\star}\leq\dfrac{\gamma_{D_k} C_k I_k(0)}{\mu_k},$$ by using the fact that $D_k(t)=\displaystyle\int_0^t \gamma_{D_k}I_k(\tau)d\tau$. Hence, by an appropriate choice of the parameters $\alpha_1^k$ and $\alpha_2^k$, it is possible to make $\dfrac{C_k}{\mu_k}$ small, which has for consequence to make $D_k^{\star}$ small.\\
Furthermore, remark that the trajectory leaves $\mathcal{D}$ if there exists $k$ such that $S_k =0$ or $I_k=0$. \textcolor{black}{This event will occur asymptotically, meaning that as time tends to infinity, it is inevitable that one of these conditions will be satisfied.} In those case, the state feedback is not well defined. Therefore, an adapted feedback law, based on \eqref{control}, is introduced in Section \ref{Sub_saturated}.\\
	\subsubsection{Non-negativity of the input}\phantom{i}\\ \\
	Obviously a control law of vaccination type should be described by a non-negative function. This in turn will ensure both a physical meaning and the non-negativity of the state, as stated in Proposition \ref{prop_pos_state}. The following theorem provides sufficient conditions for the non-negativity of the control law.\\
	\begin{thm}\label{thm_pos_unsaturated_fdb}
		Define  \begin{equation}\Gamma =\underset{k=1,...,n}{\max} \left(\gamma_{R_k}+\gamma_{D_k}\right)\label{Gamma_maj}\end{equation}
		For all $k=1,...,n,$ select $\alpha_1^k$ and $\alpha_2^k$ such that
		\begin{align}
			\alpha_1^k&> \left(\gamma_{R_k}+\gamma_{D_k}\right)\left(\Gamma+\displaystyle\sum_{j=1}^n  M_{kj}\right) \label{alpha1}\\
			\intertext{and}
			\alpha_2^k&=\gamma_{R_k}+\gamma_{D_k}+\Gamma+\displaystyle\sum_{j=1}^n M_{kj}.\label{alpha2}
		\end{align}	
		Then\textcolor{black}{, as long as the closed-loop states with the state feedback \eqref{control} remain in $\mathcal{D},$} the set $B$ is forward invariant for model \eqref{SIRD} with input \eqref{control}, the input generated by the control law \eqref{control} is non-negative andthe conclusions of Theorem \ref{thm_stab_exp} hold.
	\end{thm}\vspace{0.3cm}
	\begin{proof}
		First, one can notice that $S_k(t)$ and $I_k(t)$ 
		are non-negative for all $t\geq 0$ and $k=1,...,n$. This follows a same reasoning as for Proposition \ref{prop_pos_state} (since the essential non-negativity of the functions $f_k$ is independent of the choice of $\theta_k$ for those variables). \\
Furthermore, to show that  $\theta_k=\dfrac{v_k-L_f^2h_k}{L_{g_k}L_fh_k}\geq 0$ it is enough to show that $\tilde{\theta}_k=-v_k+L_f^2h_k\geq 0$ since $L_{g_k}L_fh_k=L_{g_k}f_k=p_kf_{n+k}< 0$ for $x\in\mathcal{D}$. Notice that, for ease of readability, the dependence in $t$ or in $x(t)$ is dropped in the following calculation. It remains to study the sign of 	$\tilde{\theta}_k.$ First observe that  
		\begin{align*}
	\tilde{\theta}_k&=\lambda S_k \displaystyle\sum_{j=1}^{n}C_{kj}f_j-\left(\gamma_{R_k}+\gamma_{D_k}\right)f_k+\lambda f_{n+k}	\displaystyle\sum_{j=1}^{n}C_{kj}I_j\\ & \hspace{0.5cm}+\alpha_2^kf_k+\alpha_1^k I_k,\\
	&=\lambda S_k \displaystyle\sum_{j=1}^{n}C_{kj}f_j+f_k\left(\Gamma+\displaystyle\sum_{j=1}^{n}M_{kj}\right)+\lambda f_{n+k}	\displaystyle\sum_{j=1}^{n}C_{kj}I_j\\&\hspace{0.5cm}+\alpha_1^k I_k,\\
	\intertext{by the choice of condition \eqref{alpha2}. The definition of $f_{n+k}$ in \eqref{f_nk} implies that}
	\tilde{\theta}_k&=\lambda S_k \displaystyle\sum_{j=1}^{n}C_{kj}f_j+f_k\left(\Gamma+\displaystyle\sum_{j=1}^{n}M_{kj}\right)-\lambda^2S_k\left(	\displaystyle\sum_{j=1}^{n}C_{kj}I_j\right)^2\\
	&\hspace{0.5cm}+\alpha_1^k I_k\\
	&\geq\lambda S_k \displaystyle\sum_{j=1}^{n}C_{kj}f_j+f_k\left(\Gamma+\displaystyle\sum_{j=1}^{n}M_{kj}\right)-\lambda^2S_k\left(	\displaystyle\sum_{j=1}^{n}C_{kj}I_j\right)^2\\
	&\hspace{0.5cm}+ \left(\gamma_{R_k}+\gamma_{D_k}\right)\left(\Gamma+\displaystyle\sum_{j=1}^n  M_{kj}\right) I_k\\ 
	\intertext{by the choice of condition \eqref{alpha1}. Using definition of $f_k$, it follows that}
\tilde{\theta}_k&\geq\lambda S_k \displaystyle\sum_{j=1}^{n}C_{kj}f_j-\lambda^2S_k\left(	\displaystyle\sum_{j=1}^{n}C_{kj}I_j\right)^2\\
&\hspace{0.5cm}+\left(\lambda S_k\displaystyle\sum_{j=1}^{n}C_{kj}I_j-\left(\gamma_{R_k}+\gamma_{D_k}\right)I_k\right)\left(\Gamma+\displaystyle\sum_{j=1}^{n}M_{kj}\right) \\ 	
&\hspace{0.5cm}+\left(\gamma_{R_k}+\gamma_{D_k}\right)\left(\Gamma+\displaystyle\sum_{j=1}^n  M_{kj}\right) I_k\\
&\geq\lambda S_k\left(\displaystyle\sum_{j=1}^{n}C_{kj}f_j+\left(\displaystyle\sum_{j=1}^{n}C_{kj}I_j\right)\left(\Gamma+\displaystyle\sum_{j=1}^{n}M_{kj}\right)\right.\\ 	
&\hspace{0.5cm}-\left.\lambda\left(	\displaystyle\sum_{j=1}^{n}C_{kj}I_j\right)^2\right)\\
&\geq\lambda S_k\left(\displaystyle\sum_{j=1}^{n}C_{kj}\left(\lambda S_j\displaystyle\sum_{l=1}^nC_{jl}I_l-\Gamma I_l\right)\right.\\&\hspace{1cm}+\left.\left(\displaystyle\sum_{j=1}^{n}C_{kj}I_j\right)\left(\Gamma+\displaystyle\sum_{j=1}^{n}M_{kj}\right)-\lambda\left(	\displaystyle\sum_{j=1}^{n}C_{kj}I_j\right)^2\right)
\intertext{thanks to the definition of $f_j$ in \eqref{f_k} and of $\Gamma$ in \eqref{Gamma_maj}. Thus}
\tilde{\theta}_k&\geq\lambda S_k\left(\displaystyle\sum_{j=1}^{n}C_{kj}\left(\lambda S_j\displaystyle\sum_{l=1}^nC_{jl}I_l\right)\right.\\&\hspace{1cm}+\left.\left(\displaystyle\sum_{j=1}^{n}C_{kj}I_j\right)\left(\displaystyle\sum_{j=1}^{n}M_{kj}\right)-\lambda\left(	\displaystyle\sum_{j=1}^{n}C_{kj}I_j\right)^2\right)\\
&=\lambda S_k\left(\displaystyle\sum_{j=1}^{n}C_{kj}\left(\lambda S_j\displaystyle\sum_{l=1}^nC_{jl}I_l\right)\right.\\&\hspace{1cm}+\left.\left(\displaystyle\sum_{j=1}^{n}C_{kj}I_j\right)\left(\displaystyle\sum_{j=1}^{n}M_{kj}-\lambda\displaystyle\sum_{j=1}^{n}C_{kj}I_j\right)\right)\\
&\geq\lambda S_k\left(\displaystyle\sum_{j=1}^{n}C_{kj}\left(\lambda S_j\displaystyle\sum_{l=1}^nC_{jl}I_l\right)\right.\\&\hspace{1cm}+\left.\left(\displaystyle\sum_{j=1}^{n}C_{kj}I_j\right)\left(\displaystyle\sum_{j=1}^{n}M_{kj}-\lambda\displaystyle\sum_{j=1}^{n}M_{kj}\right)\right),\\
\intertext{since $C_{kj}I_j\leq C_{kj}N_j=M_{kj}$. Finally,}
\tilde{\theta}_k&\geq\lambda S_k\left(\displaystyle\sum_{j=1}^{n}C_{kj}\left(\lambda S_j\displaystyle\sum_{l=1}^nC_{jl}I_l\right)\right.\\&\hspace{1cm}+\left.\left(1- \lambda\right)\displaystyle\sum_{j=1}^{n}C_{kj}I_j\displaystyle\sum_{j=1}^{n}M_{kj}\right)\geq 0
\end{align*}
since $\lambda\leq 1$.\\
Moreover, since the feedback design parameters are positive, Theorem \ref{thm_stab_exp} concludes the proof.\end{proof}	

	\subsection{Constrained Feedback law}\label{Sub_saturated}
\begin{figure}
		\centering
		\includegraphics[scale=0.3]{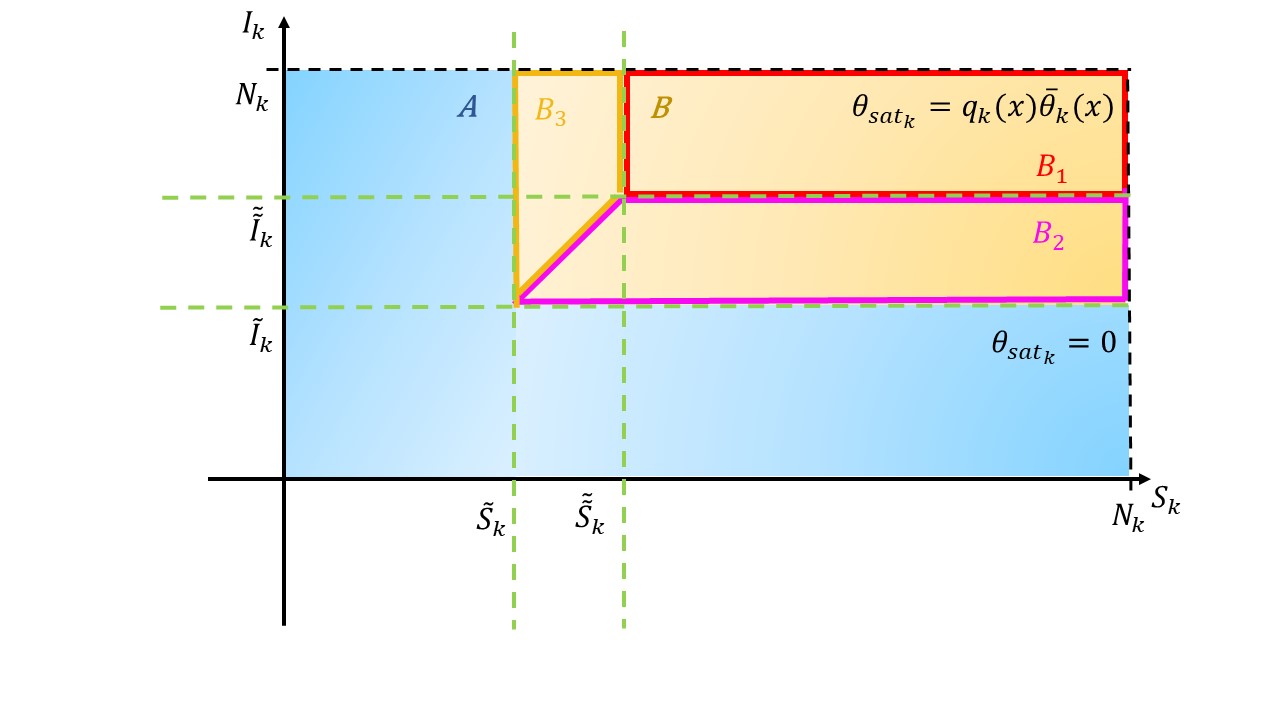}
		\caption{Representation of the construction of the constrained feedback-law, $\theta_{sat_k}(x)$}  
		\label{fig_lipschitz} 
	\end{figure}

\begin{figure}
		\centering
		\includegraphics[scale=0.3]{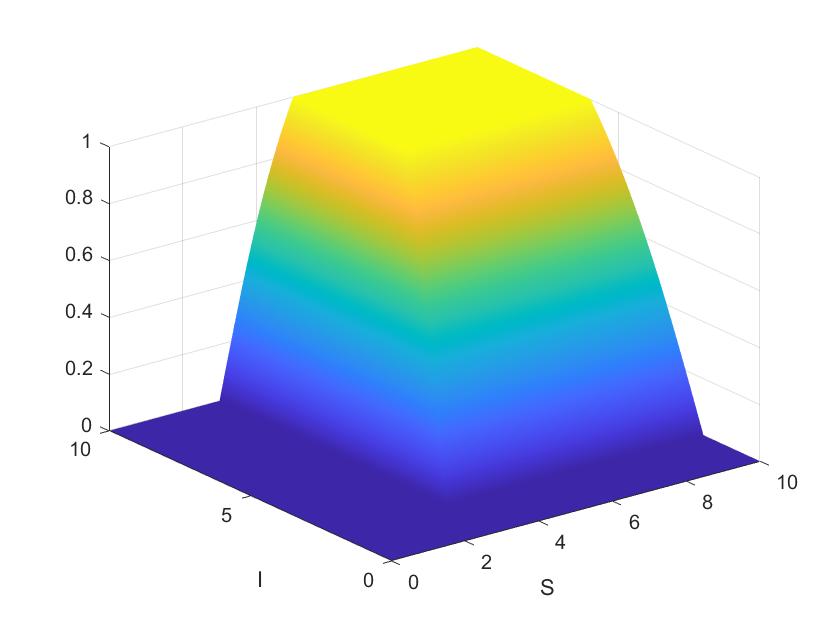}
		\caption{Representation of the control law  $\theta_{sat_k}(t)$}  
		\label{fig_controle_sat} 
	\end{figure}

	In the previous section a stabilizing state feedback law has been defined by $$\theta\left(t\right)=\begin{pmatrix}
		\dfrac{v_1\left(x\left(t\right)\right)-L_f^2h_1\left(x\left(t\right)\right)}{L_{g_1}L_fh_1\left(x\left(t\right)\right)} &
		\cdots &
		\dfrac{v_n\left(x\left(t\right)\right)-L_f^2h_n\left(x\left(t\right)\right)}{L_{g_n}L_fh_n\left(x\left(t\right)\right)}
	\end{pmatrix}^T,$$ where $L_{g_k}L_fh_k\left(x\left(t\right)\right)=-p_k\lambda S_k\left(t\right)\displaystyle\sum_{j=1}^n C_{kj}I_j\left(t\right).$
	As predicted in Theorem \ref{thm_stab_exp}, $\displaystyle\sum_{j=1}^n C_{kj}I_j\left(t\right)$ tends to $0$ and the feedback blows up. Moreover, there is no proof that the solution will remain in $\mathcal{D}$, which is an assumption needed to obtain the conclusion of Theorem~\ref{thm_pos_unsaturated_fdb}. To avoid this and to take into account design specifications, where an amplitude constraint on the control, denoted by $\theta_{sup}$, is imposed in practice, we define a new control law inspired by the previous one but with saturation and insurance that $x$ is in $\mathcal{D}$ when the law \eqref{control} is used. This law will solve the second problem that consists, as mentioned before, of designing an amplitude limited control that improves performance with respect to the open-loop system regarding the peak of total infected individuals (thanks to Proposition \ref{prop_goal_vacc}, at least in the case of one class of age), while maintaining asymptotic convergence. \\ \\
	Let us define a new state feedback law by $$	\theta_{sat_k} : \left[0,N_k\right]^{n+1} \to \mathbb{R}$$ where, with $x=(S_k,I_1...I_n)$, \begin{equation}\label{control_sat}
		\theta_{sat_k}(x) = \left\{
		\begin{array}{ll}
			(q_k\bar{\theta}_k)(x) & \mbox{if } S_k\geq \tilde{S}_k \text{ and } I_k\geq \tilde{I}_k  \text{ ($B$ area)} \vspace{0.3cm}\\
			0 & \mbox{otherwise } \text{ ($A$ area)}
		\end{array}
		\right.
	\end{equation}
	where $(q_k\bar{\theta}_k)(x)$ denotes $q_k(x)\bar{\theta}_k(x)$ and the area refers to the visual representation, available in Figure \ref{fig_lipschitz}.
	The functions $\bar{\theta}_k$ and $q_k$ are defined on $\left[\tilde{S}_k,N_k\right]\times\left[0,N_k\right]\times \cdots \times \left[0,N_k\right] \times \left[\tilde{I}_k,N_k\right]\times \left[0,N_k\right]\times \cdots \times \left[0,N_k\right]$ and allow the new state feedback law \eqref{control_sat} to have desired properties, such as boundedness and Lipschitz properties\textcolor{black}{, which are proved in the sequel.} They are given by
	\begin{equation}\label{fct_bar_theta_k}
		\bar{\theta}_k(x)= \left\{
		\begin{array}{ll}
			\theta_k(x) & \mbox{if } \theta_k(x)\leq \theta_{sup} \vspace{0.3cm}\\
			\theta_{sup} &  \mbox{otherwise, }\\
		\end{array}
		\right.
	\end{equation}
where $\theta_{sup}$ correspond to an amplitude constraint on the control and 
	\begin{equation}\label{fct_gk}
		q_k(x)=\left\{
		\begin{array}{ll}
			1 & \mbox{if } S_k\geq \tilde{\tilde{S}}_k \text{ and } I_k\geq \tilde{\tilde{I}}_k\\ & \text{($B_1$ area)} \vspace{0.3cm}\\
			\dfrac{4}{\pi}\arctan\left(\dfrac{I_k-\tilde{I}_k}{\tilde{\tilde{I}}_k-\tilde{I}_k}\right) & \mbox{if }I_k\leq \\ &\dfrac{ \tilde{\tilde{I}}_k-\tilde{I}_k}{\tilde{\tilde{S}}_k-\tilde{S}_k}\left(S_k-\tilde{S}_k\right)+\tilde{I}_k\\ &\text{and } I_k<\tilde{\tilde{I}}_k\text{ ($B_2$ area)}\vspace{0.3cm}\\
			\dfrac{4}{\pi}\arctan\left(\dfrac{S_k-\tilde{S}_k}{\tilde{\tilde{S}}_k-\tilde{S}_k}\right) & \mbox{otherwise}\text{ ($B_3$ area)},\\
		\end{array}
		\right.
	\end{equation}
	where the constants $\tilde{\tilde{S}}_k$ and $\tilde{\tilde{I}}_k$ are chosen larger than $\tilde{S}_k$ and $\tilde{I}_k$, respectively. Moreover, $\tilde{S}_k$ and $\tilde{I}_k$  have to be selected appropriately as shown below.\\
Note that this law is constructed following a similar approach to the one introduced in \cite{Molnar_2021}, where the authors also propose a constrained control law but to ensure the forward invariance of a safe set. This safe set can represent the fact that the number of infected individuals remains below a certain threshold.\\
	This new law is represented in Figure \ref{fig_controle_sat} for $\tilde{S}=3,	\tilde{\tilde{S}}=5,
	\tilde{I}=2,
	\tilde{\tilde{I}}=6$. With this definition, $\theta_{sat_k}$ is globally bounded and Lipschitz as \textcolor{black}{a} function of $I_k$ and $S_k$. Those properties will be useful in the design of the observer-based output feedback, developed in Section \ref{Sect_output_fdb}. Moreover, $\theta_{sat_k}$ has the advantage to be well defined even if $x$ does not belong to $\mathcal{D}$ anymore. Indeed, it is based on the ``old'' law \eqref{control}, denoted by $\theta$, when there is the absolute certainty that $x$ belongs to $\mathcal{D}$ and is zero otherwise, as defined in \eqref{control_sat}. Some performance properties of this law are shown in the sequel.
	
	The first result states the asymptotic stability of the trajectories for the closed-loop model under saturated feedback.\\
	\begin{prop}\label{thm_stab_sat_vacc}
	Choose the control tuning parameters to satisfy the assumptions of Theorem \ref{thm_pos_unsaturated_fdb}. Then, the saturated state feedback \eqref{control_sat} implies the pointwise attractivity of the disease-free equilibria set of model \eqref{SIRD}. Moreover the susceptible, recovered and deceased populations converge asymptotically to some constants $S_k^{\star}, R_k^{\star}$ and $D_k^{\star}$, respectively, for $k=1,...,n$.
	\end{prop}\vspace{0.3cm}
	\begin{proof}
		This proposition is a direct consequence of Proposition \ref{prop_asympt_stab} since the saturated state feedback \eqref{control_sat} is, by construction, non-negative.
	\end{proof}\vspace{0.3cm}
	The following lemma introduces an invariant set helpful to prove that the constrained feedback law \eqref{control_sat} has a finite number of ``jumps" \textcolor{black}{in the definition \eqref{control_sat} of $\theta_{sat_k}$. This implies that there exists a time $T$ such that $x(t)$ remains in the $A$ area or in the $B$ area, for $t\geq T$}.
	\begin{lemma}\label{lemma_invariant}
		Let $\tilde{I}_k>0$ and $\tilde{S}_k\leq \dfrac{\left(\gamma_{R_k}+\gamma_{D_k}\right)\tilde{I}_k}{\lambda\displaystyle\sum_{j=1}^n M_{kj}}.$ 
		Then, the set $\left[0,\tilde{S}_k\right] \times \left[0,\tilde{I}_k\right]\times [0,N_k] \times [0,N_k]$ is forward invariant for model \eqref{SIRD}, for any non-negative input.
	\end{lemma}\vspace{0.3cm}
	\begin{proof}
		First, observe that the non-negativity of the states follows from Proposition \ref{prop_pos_state} since $\theta_{sat_k}$ is non-negative. The same proposition implies that each state is bounded by $N_k$. Moreover, let $(S_k(t), I_k(t))\in \left[0,\tilde{S}_k\right] \times \left[0,\tilde{I}_k\right]$, for some $t\geq 0$ arbitrarily fixed.\\
		For, $S_k$, two cases can happen. Either there exists $t_1$ such that $S_k(t_1)=0$, in that case, $\dot{S_k}(t_1)=0$. It follows that $S_k(t)=0$ for all $t\geq t_1$ and remains in $\left[0,\tilde{S}_k\right]$. Either $S_k(t)\neq 0$ for all $t\geq 0$. It follows that $\dfrac{dS_k(t)}{dt}<0$ since the input is non-negative. Therefore, at the border, when $S_k(t)=\tilde{S}_k$, the state trajectories for the susceptible remains in the set $ \left[0,\tilde{S}_k\right].$
		Furthemore, focus on what happens at the other border, when $I_k(t)=\tilde{I_k}$ and $S_k(t)\leq \tilde{S}_k(t)$. Then, \begin{align*}
			&S_k(t)\leq \tilde{S}_k\leq \dfrac{\left(\gamma_{R_k}+\gamma_{D_k}\right)\tilde{I}_k}{\lambda\sum_{j=1}^n M_{kj}}\leq \dfrac{\left(\gamma_{R_k}+\gamma_{D_k}\right)\tilde{I}_k}{\lambda\sum_{j=1}^n C_{kj}I_j}\\
			\Rightarrow & \lambda S_k(t) \displaystyle\sum_{j=1}^n C_{kj}I_j-\left(\gamma_{R_k}+\gamma_{D_k}\right)\tilde{I}_k\leq 0\\
			\Leftrightarrow & \dfrac{dI_k(t)}{dt}\leq 0.
		\end{align*}
		Therefore, we can conclude that the set $\left[0,\tilde{S}_k\right] \times \left[0,\tilde{I}_k\right]\times [0,N_k] \times [0,N_k]$ is forward invariant.
	\end{proof}\\ \\
	The next property ensures that the solution will end up in \textcolor{black}{the invariant set introduced in Lemma \ref{lemma_invariant}. Consequently, the control law \eqref{control_sat} cannot switch anymore, in finite time. In view of the application, this guarantees that there exists a finite time after which further immunization is unnecessary.} 
	\begin{property}\label{lemma_finite}
		Let $\tilde{I}_k$ and $\tilde{S}_k$ be chosen as in Lemma \ref{lemma_invariant} and let $x_0\in B$. Then, there exists a finite time $t_{x_0}$ such that $I_k(t_{x_0})=\tilde{I}_k$ and $\forall$ $t>t_{x_0}, I_k(t)<\tilde{I}_k$. Therefore, there is a finite number of switches, m, for the function $\theta_{sat_k}$ defined in \eqref{control_sat}, at times $t_1<...<t_m\leq t_{x_0}$.
	\end{property}
	\begin{figure}
		\centering
		\includegraphics[scale=0.28]{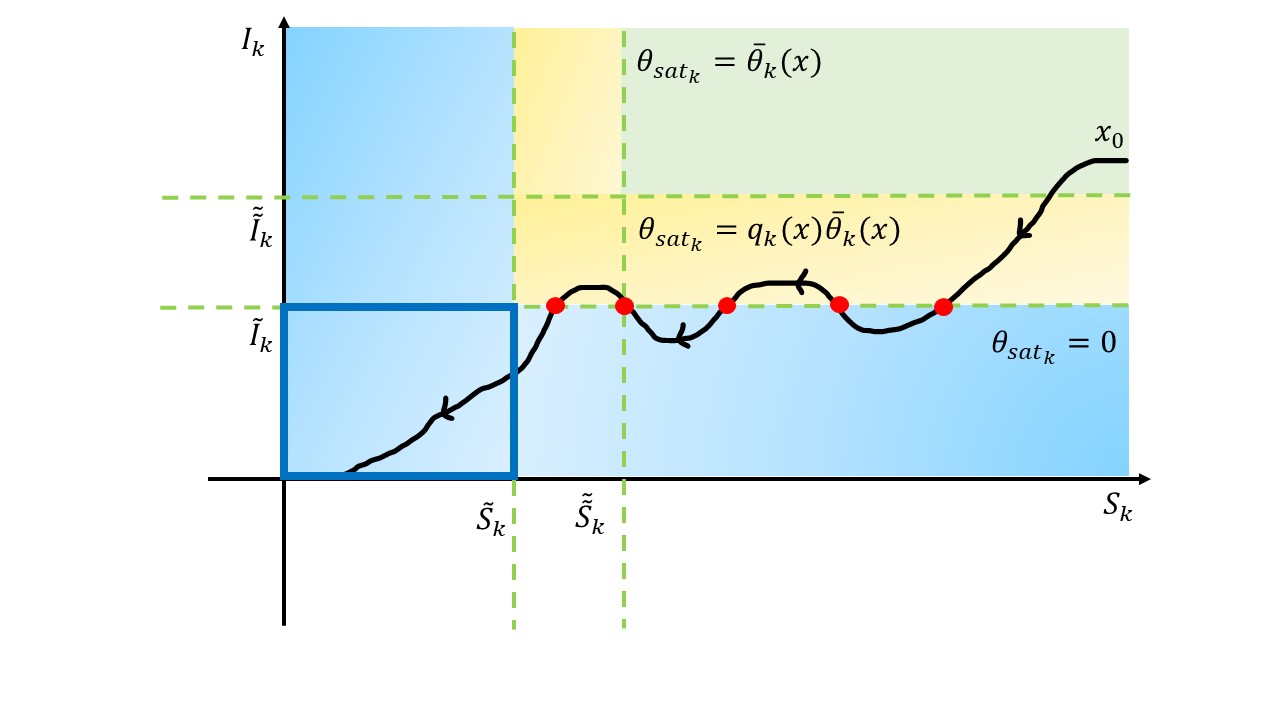}
		\caption{Illustration of the finite number of jumps for the saturated control law, $\theta_{sat_k}(t)$}  
		\label{fig_vacc_switch} 
	\end{figure}\vspace{0.3cm}
	\begin{proof}
		Let $k=1,...,n$, and consider three cases. First,  assume that there exists $T>0$ such that $S_k(T)= 0$. In that case, $S_k(t)=0$ for all $t\geq T$ since $\dot{S}_k(T)=0$, implying that $S_k(t)=0<\tilde{S}_k$ for all $t\geq T$. Then, no more jump is possible thanks to the definition \eqref{control_sat} of $\theta_{sat_k}$. \\
		Second, assume  that there exists $T>0$ such that $\lambda\displaystyle\sum_{j=1}^{n} C_{kj}I_j(T)+p_k\theta_{sat_k}(T)=0$. Thanks to the non-negativity of the elements in the equation, this case is only possible if $I_k(T)=0$ for all $k=1,...,n$ and if $\theta_{sat_k}(T)=0$. Therefore, $\dot{I}_k(T)=0$ implying that $I_k(t)=0<\tilde{I}_k$ for all $t\geq T$. Then, no more jump is possible thanks to the definition \eqref{control_sat} of $\theta_{sat_k}$.\\
		Finally, in the case where $S_k(t)\neq 0$ and $\lambda\displaystyle\sum_{j=1}^{n} C_{kj}I_j(t)+p_k\theta_{sat_k}(t)\neq 0$, $k=1,...,n$, $t\geq 0$ then $S_k(t)$ is strictly decreasing, since the feedback $\theta_{sat_k}(t)$ is non-negative. Moreover, $S_k$ has zero as a lower bound therefore it converges to some equilibrium $S_k^{\star}$. Since $I_k(t)$ tends to zero, thanks to Proposition \ref{thm_stab_sat_vacc}, so the couple $(S_k(t),I_k(t))$ tends to $(S_k^{\star},0)$ when $t$ goes to infinity. By absurd, assume that there exists an infinite number of switches. In this case, there always exist time instants $t$ such that $I_k(t)\geq \tilde{I}_k$, for the switches to occur. This contradicts the fact that the infected individuals tend to zero. Remark that, in the case where the equilibrium is such that $S_k^{\star}<\tilde{S}_k$, illustrated in Figure \ref{fig_vacc_switch}, the strictly decreasing property of $S_k$ can be invoked to show that there exists a finite time instant $T$ where $S_k(T)= \tilde{S}_k$. Afterwards, the control law will remain switched off and the trajectory will enter the invariant region $\left[0,\tilde{S}_k\right] \times\left[0,\tilde{I}_k\right]$ since the infected individuals tend to $0$.
	\end{proof}\\ \\
	Finally, it can be shown that the constrained feedback law, as defined in \eqref{control_sat}, is Lipschitz. This property is a crucial assumption for the design of an observer-based output feedback, presented in the next section.\\ 
	
	\begin{property}\label{property_Lipschitz}
		The function $\theta_{sat_k}$ is Lipschitz on its domain\textcolor{black}{, bounded and such that $\theta_{sat_k}(0)=0$.}	
	\end{property}\vspace{0.3cm}
	
	\begin{proof}
		This proof is quite computational. Therefore, only some key elements are presented here. Some additional information is available in the Appendix. Moreover, the reader can refer to Figure \ref{fig_lipschitz} to better understand the areas considered in each case.\\
		Let $z_k=(x_k,y_1,\cdots,y_k,\cdots,y_n)^T$ and $z_k'=(x_k',y'_1,\cdots,y'_k,\cdots,y'_n)^T \in \left[0,N_k\right]^{n+1}.$ Three main cases can be identified.
		\begin{enumerate}
			\item Assume that $z_k$ and $z_k'$ are such that $x_k, x_k'<\tilde{S}_k$ or $y_k,y_k'<\tilde{I}_k$ ($z_k$ and $z_k' \in A$). Hence, $\vert \theta_{sat_k}(z_k)-\theta_{sat_k}(z'_k) \vert=0\leq L_0\Vert z_k-z'_k \Vert$  for all $L_0>0$.\\
			\item Assume that $z_k$ and $z_k'$ are such that $(x_k,y_k)$ and $(x_k',y_k')\in\left[\tilde{S}_k,N_k\right]\times \left[\tilde{I}_k,N_k\right]$ ($z_k$ and $z_k' \in B$). Therefore, \begin{align*}
				\vert \theta_{sat_k}(z_k)-\theta_{sat_k}(z'_k) \vert&= \vert (q_k\bar{\theta}_k)(z_k)-(q_k\bar{\theta}_k)(z'_k) \vert\\ &\leq L\Vert z_k-z'_k \Vert,
			\end{align*} with $L>0$, the Lipschitz constant of $q_k\bar{\theta}_k$ (see details in the Appendix).\\
			\item Assume that $z_k$ and $z_k'$ are such that $(x_k,y_k)\in\left[\tilde{S}_k,N_k\right]\times \left[\tilde{I}_k,N_k\right]$ and $x_k'<\tilde{S}_k$ or $y_k'<\tilde{I}_k$ ($z_k\in B$ and $z_k' \in A$). Then,
			\begin{align*}
				\vert \theta_{sat_k}(z_k)-\theta_{sat_k}(z'_k) \vert&= \vert (q_k\bar{\theta}_k)(z_k) \vert\\ &=(q_k\bar{\theta}_k)(z_k).
			\end{align*}
			Three subcases can again be identified.
			\begin{enumerate}
				\item $z_k$ such that $(x_k,y_k)\in\left[\tilde{S}_k,N_k\right]\times \left[\tilde{I}_k,N_k\right]$ with $x_k\geq \tilde{\tilde{S}}_k$ and $y_k\geq \tilde{\tilde{I}}_k$ ($z_k \in B_1$). Hence,
				\begin{align*}
					\vert \theta_{sat_k}(z_k)-\theta_{sat_k}(z'_k) \vert&=\bar{\theta}_k(z_k)\\
					&\leq \theta_{sup}\\
					&\leq L_1\Vert z_k-z_k'\Vert
				\end{align*}
				where $L_1=\dfrac{\theta_{sup}}{d(A,B_1)}$ with $$d(A,B_1)=\underset{a\in A, b\in B_1}{\inf}\Vert a-b\Vert,$$ since by definition, $d(A,B_1)\leq \Vert z_k-z'_k\Vert$. 
				\item $z_k$ such that $(x_k,y_k)\in\left[\tilde{S}_k,N_k\right]\times \left[\tilde{I}_k,N_k\right]$ with $y_k\leq \dfrac{ \tilde{\tilde{I}}_k-\tilde{I}_k}{\tilde{\tilde{S}}_k-\tilde{S}_k}\left(x_k-\tilde{S}_k\right)+\tilde{I}_k $ and $y_k<\tilde{\tilde{I}}_k$ ($z_k \in B_2$).	
				Therefore,
				\begin{align*}
					\vert \theta_{sat_k}(z_k)-\theta_{sat_k}&(z'_k) \vert=q_k(z_k)\bar{\theta}_k(z_k)\\
					&\leq \dfrac{4}{\pi}\arctan\left(\dfrac{y_k-\tilde{I}_k}{\tilde{\tilde{I}}_k-\tilde{I}_k}\right)\theta_{sup}\\
					&\leq \dfrac{4 \theta_{sup}}{\pi(\tilde{\tilde{I}}_k-\tilde{I}_k)}(y_k-\tilde{I}_k)\\
					&=: L_2 (y_k-\tilde{I}_k)
				\end{align*}
				since $\arctan$ is a Lipschitz function with constant $1$. Now, if $y_k'\leq \tilde{I_k}$, then
				\begin{align*}
					\vert \theta_{sat_k}(z_k)-\theta_{sat_k}(z'_k) \vert
					&\leq L_2\left[(y_k-y_k')+(y_k'-\tilde{I}_k)\right]\\
					&\leq L_2(y_k-y_k')\\
					&\leq L_2\Vert z_k-z'_k\Vert.
				\end{align*}
				On the other hand, if $y_k'> \tilde{I_k}$, that means that $x_k'\leq \tilde{S}_k$ and we have that 
				\begin{align*}
					\vert \theta_{sat_k}(z_k)-\theta_{sat_k}&(z'_k) \vert
					\leq L_2\dfrac{ \tilde{\tilde{I}}_k-\tilde{I}_k}{\tilde{\tilde{S}}_k-\tilde{S}_k}(x_k-\tilde{S}_k)\\
					&=:L_3\left[(x_k-x_k')+(x_k'-\tilde{S}_k)\right]\\
					&\leq L_3(x_k-x_k')\\
					&\leq L_3\Vert z_k-z'_k\Vert.
				\end{align*}
				\item $z_k$ such that $(x_k,y_k)\in\left[\tilde{S}_k,N_k\right]\times \left[\tilde{I}_k,N_k\right]$ with $y_k> \dfrac{ \tilde{\tilde{I}}_k-\tilde{I}_k}{\tilde{\tilde{S}}_k-\tilde{S}_k}\left(x_k-\tilde{S}_k\right)+\tilde{I}_k $ and $x_k<\tilde{\tilde{S}}_k$ ($z_k \in B_3$). By using similar reasoning as (b), it follows that  $$\vert \theta_{sat_k}(z_k)-\theta_{sat_k}(z'_k) \vert\leq L_4 \Vert z_k-z'_k\Vert. $$	
			\end{enumerate}
		\end{enumerate}
		Therefore, there exists $K=\max\left\{L_0,L,L_1,L_2,L_3,L_4\right\}>0$ such that, for all  $z_k$ and $z_k'\in \left[0,N_k\right]^{n+1},$ $$  \vert \theta_{sat_k}(z_k)-\theta_{sat_k}(z'_k) \vert\leq K \Vert z_k-z'_k\Vert.$$
	\textcolor{black}{Moreover, $\theta_{sat_k}(0)=0$ by definition since $\theta_{sat_k}(0)$ corresponds to the case where $S_k=I_k=0$. The boundedness follows from the saturation imposed in the definition of $\theta_{sat_k}$.  }	
	\end{proof}

	\section{Observer-based output feedback design}\label{Sect_output_fdb}
	\label{Observer}
	\textcolor{black}{The feedback considered in the previous section may need the use of state variables which are not measured. Typically indeed, in epidemics, one usually measures} the number of  dead individuals in each class of age. Therefore, the goal of this section is to design and analyze an observer-based output feedback law which is based only on such measurements. The technique used in this work is inspired by the separation principle introduced in \cite{atassi_1999} for multi-input multi-output systems. \\ \\
First, a change of variable is introduced to rewrite model \eqref{SIRD} as a chain of integrators. Furthermore, as mentioned previously, in the case of epidemic, data on the dead individuals can be collected. Therefore, in the following, we use model \eqref{SIRD} where only the equations of $S_k, I_k$ and $D_k$ are kept due to the constant population assumption, with outputs given by $y_k(t)=D_k(t)$, $k=1,...,n$. Hence, the following change of variable is performed: 
	\begin{align}\label{change_var}
		z_{1k}&=D_k,\nonumber\\
		z_{2k}&=\gamma_{D_k}I_k\nonumber\\
		z_{3k}&=\gamma_{D_k}\left\{\lambda S_k\left(t\right) \displaystyle\sum_{j=1}^n C_{kj}I_j(t) -\left(\gamma_{R_k}+\gamma_{D_k}\right)I_k\left(t\right)\right\}
	\end{align}
	
	\noindent for $k=1,...,n$.\\ \textcolor{black}{Letting $z=\begin{pmatrix}z_{11}& z_{21}& z_{31}& z_{12}& z_{22}& z_{32}& ...& z_{3n}\end{pmatrix}^T$, this change of variables implies that $z\in \mathcal{Z}$ where \begin{align}\label{mathcal_Z}
		\mathcal{Z}:=&\biggl\{z\in \mathbb{R}^{3n}: 0\leq z_{1k}\leq N_k,0\leq z_{2k}\leq \gamma_{D_k}N_k,\nonumber\\ &-\left(\gamma_{R_k}+\gamma_{D_k}\right)z_{2k} \leq z_{3k}\\ &\leq N_k\gamma_{D_k}\lambda\displaystyle\sum_{j=1}^n \tilde{C}_{kj}z_{2j}-\left(\gamma_{R_k}+\gamma_{D_k}\right)z_{2k} \hspace{0.2cm} \forall k=1,...,n\biggr\},\nonumber
		\end{align} corresponds to the set $B$ for the original variables.} Remark that this change of variable corresponds to the change of variable \eqref{change_of_variable} introduced previously, for $z_{2k}$ and $z_{3k}$, $k=1,...,n$, up to multiplication by the factor $\gamma_{D_k}.$ Therefore, it is also invertible if $\displaystyle\sum_{j=1}^n C_{kj}I_j\neq 0$. Indeed, the nominal coordinates can be obtained by 
	\begin{align}\label{inverse_change_var}
		D_k&=z_{1k},\nonumber\\
		I_k&=\dfrac{z_{2k}}{\gamma_{D_k}}\\
		S_k&= \left\{\begin{array}{ll}
			S_k^{\star} & \text{ if }z_{2j}=0 \hspace{0.2cm}\forall j\vspace{0.3cm}\\
			\dfrac{\dfrac{z_{3k}}{\gamma_{D_k}}+\left(\gamma_{R_k}+\gamma_{D_k}\right)\dfrac{z_{2k}}{\gamma_{D_k}}}{\lambda\displaystyle\sum_{j=1}^n \tilde{C}_{kj} z_{2j}} & \text{otherwise}
			\nonumber
		\end{array}\right.
	\end{align}
	where $\tilde{C}_{kj}=\dfrac{C_{kj}}{\gamma_{D_j}}=\dfrac{M_{kj}}{\gamma_{D_j}N_j}$ and $S_k^{\star}$ is the equilibrium for the susceptible individuals in the kth class of age. Notice that this value does not need to be computed since it is known that the state-feedback is $0$ in this case. In simulation, one does not need to compute the output feedback in this case, it can just be set to $0$.\\
	One can notice that $z_{2j}=0$ for all $j=1,...,n$ is equivalent to the fact that $I_j= 0$ for all $j=1,...,n$. Therefore, it implies that $z_{3j}=0$ for all $j=1,...,n$ in view of the change of variable \eqref{change_var}.
	
	Therefore, in the new variables, model \eqref{SIRD} becomes
	\begin{subequations}\label{SIRD_z}
		\begin{equation}\label{SIRD_za}
			\left\{\begin{array}{rl}
				\dot{z}_{1_k}\left(t\right) &= z_{2k}\left(t\right) \\ 
				\dot{z}_{2_k}\left(t\right) &=z_{3k}\left(t\right)\vspace{0.2cm}\\ 
				\dot{z}_{3_k}\left(t\right)&=\left\{\begin{array}{ll}
					0 & \text{if } z_{2j}, z_{3j}=0\hspace{0.2cm} \forall j\vspace{0.2cm}\\
					h_k(z(t),\textcolor{black}{u}(z(t))) & \text{otherwise}
				\end{array}\right.
			\end{array}
			\right.	
		\end{equation}
		where
		\begin{align}
			\begin{split}
				h_k(z(t),\textcolor{black}{u}(z(t)))=	\left(z_{3k}(t)+\left(\gamma_{R_k}+\gamma_{D_k}\right)z_{2k}(t)\right)\phantom{iiiiiiiiiiiiiiiiiiiiiii}\\  \times \left[\lambda \displaystyle\sum_{j=1}^n \tilde{C}_{kj}z_{2j}(t)-p_k\textcolor{black}{u}_k(z(t))+\dfrac{\displaystyle\sum_{j=1}^n \tilde{C}_{kj}z_{3j}(t)}{\displaystyle\sum_{j=1}^n \tilde{C}_{kj}z_{2j}(t)}-\left(\gamma_{R_k}+\gamma_{D_k}\right)\right]\phantom{iiiiii}\\+\left(\gamma_{R_k}+\gamma_{D_k}\right)^2z_{2k}(t)\phantom{iiiiiiiiiiiiiiiiiiiiiiiiiiiiiiiiiiiiiiiiiiiiiiiiiiiiiiiiiiiii}
			\end{split}
		\end{align}
	\end{subequations}
\textcolor{black}{where $u(z)=\left(u_1(z),...,u_n(z)\right)^T$ with, 	\begin{equation*}
		u_k(z)= \left\{
		\begin{array}{ll}
			\bar{\theta}_{sat_k}(z) & \mbox{if } z\in \mathcal{Z}, \vspace{0.3cm}\\
			0 &  \mbox{otherwise}\\
		\end{array}
		\right.
\end{equation*}
and $\bar{\theta}_{sat_k}$ is obtained from the $kth$ component of the given bounded non-negative state feedback \eqref{control_sat}, combined with the change of variables \eqref{inverse_change_var}.\\
This control is introduced to be defined even if the observer states leave $\mathcal{Z}$. 
\begin{rem}\label{rem_lip}
Thanks to Property \ref{property_Lipschitz} and its definition, $u_k$ is Lipschitz and bounded.
\end{rem}}
\vspace{0.2cm}
Moreover, it can be shown that the vector field of this system is continuous. Indeed, in the proof of Lemma \ref{lemma_ass_1}, it is shown that 	$\displaystyle\left\lvert\frac{\displaystyle\sum_{j=1}^n \tilde{C}_{kj}z_{3j}}{\displaystyle\sum_{j=1}^n C_{kj}z_{2j}}\displaystyle\right\rvert\leq K$. The squeeze theorem completes the argumentation.\\
Model \eqref{SIRD_z} rewrites as
		\begin{equation}\label{state-feed}
		\left\{\begin{array}{rl}
			\dot{z}\left(t\right) &= Az\left(t\right)+B\phi\left(z\left(t\right),u\left(t\right)\right) \\ 
			y\left(t\right) &=Hz\left(t\right)
		\end{array}\\
		\right.
	\end{equation}	
	with $z(0)=z_0$, where $z\in \mathcal{Z}\subseteq\mathbb{R}^{3n}$ is the state vector, \begin{align}\label{state_fdb}
		u=\theta(z(t)),
	\end{align}such that $u\in \mathcal{U}\subseteq\mathbb{R}^n,$ is the control input (i.e. the exact state feedback) and $y\in \mathcal{Y}\subseteq\mathbb{R}^n$ is the measured output. The matrices $A$, $B$ and $H$ are given by $A=\text{blockdiag}[\tilde{A},...,\tilde{A}]_{3n\times 3n}$, $B=\text{blockdiag}[\tilde{B},...,\tilde{B}]_{3n\times n}$ and $H=\text{blockdiag}[\tilde{H},...,\tilde{H}]_{n\times 3n}$, where $$\tilde{A}=\begin{pmatrix}
		0&1&0\\0&0&1\\0&0&0
	\end{pmatrix}, \tilde{B}=\begin{pmatrix}
		0&0&1
	\end{pmatrix}^T, \tilde{H}=\begin{pmatrix}
		1&0&0
	\end{pmatrix}.$$\vspace{0.3cm}	

In order to establish the main result of this section, namely Theorem \ref{thm_HG}, some additional properties are needed. They are stated in the following lemmas. Then, the desired observer-based output feedback will be obtained using the high-gain observer \eqref{HG}. For this proof the uniform norm on $\mathbb{R}^{3n}$ will be used.\vspace{0.3cm}
	\begin{lemma}\label{lemma_ass_1}
		\textcolor{black}{Let $$\mathcal{U}=\left\{U=\left(U_1,...,U_n\right)^T\in\mathbb{R}^ n : 0\leq U_k\leq \theta_{sup}, k=1,...,n\right\}.$$} The function $\phi:\mathcal{Z}\times\mathcal{U} \to \mathbb{R}^n$ is Lipschitz in its arguments on its domain. 
	\end{lemma}\vspace{0.3cm}
	\begin{proof}
		First, observe that $\forall k= 1,...,n$ and $\forall z$ such that $z_{2j}\neq 0$ and $z_{3j}\neq 0$  for some $j\in\left\{1,...,n\right\}$, $\phi_k(x,u)$ is locally Lipschitz in $z$ and in $u$ on the set
	$\mathcal{Z}$, defined in \eqref{mathcal_Z}, and $\mathcal{U}$ since $\phi_k$ is of class $C^1$ on $\left\{z\in\mathcal{Z}:z_{2j}\neq 0 \text{ or } z_{3j}\neq 0 \text{ for some }j\in\left\{1,...,n\right\}\right\}$ and on $\mathcal{U}$.\\
		Moreover, for all $z_0=\begin{pmatrix}
			z_{0_{11}} &\cdots& z_{0_{1n}} & 0 &\cdots &0
		\end{pmatrix}^T$, $\phi_k(z_0,U)$ is Lipschitz in $\mathcal{Z}$. This means that there exists $M>0$ such that $\forall z, z_0\in\mathcal{Z}$, $\vert \phi_k(z,U)-\phi_k(z_0,U)\vert \leq M \Vert z-z_0\Vert.$ Indeed, for $z\in\mathcal{Z}$ such that $z$ equals to $z_0$, the result is trivial. In addition, considering $z\neq z_0$, it follows from the definition of $\phi_k$ that \begin{align}
			\vert \phi_k(z,U)-\phi_k(z_0,U)\vert = \vert \phi_k(z,U)\vert,\phantom{iiiiiiiiiiiiiiiiiiiiiiiii}\nonumber\\
			\leq\vert\left(z_{3k}+\left(\gamma_{R_k}+\gamma_{D_k}\right)z_{2k}\right)\vert\left[\lambda \displaystyle\sum_{j=1}^n \tilde{C}_{kj}\vert z_{2j}\vert+\vert p_k U_k\vert\phantom{iii} \right.\nonumber\\\left.+\displaystyle\left\lvert\dfrac{\displaystyle\sum_{j=1}^n \tilde{C}_{kj}z_{3j}}{\displaystyle\sum_{j=1}^n \tilde{C}_{kj}z_{2j}}\displaystyle\right\rvert+\left(\gamma_{R_k}+\gamma_{D_k}\right)\right]+\left(\gamma_{R_k}+\gamma_{D_k}\right)^2\vert z_{2k}\vert,\nonumber\\
		\end{align}
		\begin{align}
			\leq\left(1+\left(\gamma_{R_k}+\gamma_{D_k}\right)\right) \Vert z-z_0\Vert \left[\lambda \displaystyle\sum_{j=1}^n \tilde{C}_{kj}\vert z_{2j}\vert+\vert p_k U_k\vert \phantom{i}\right.\nonumber\\\left.+\displaystyle\left\lvert\dfrac{\displaystyle\sum_{j=1}^n \tilde{C}_{kj}z_{3j}}{\displaystyle\sum_{j=1}^n \tilde{C}_{kj}z_{2j}}\displaystyle\right\rvert+\left(\gamma_{R_k}+\gamma_{D_k}\right)\right]+\left(\gamma_{R_k}+\gamma_{D_k}\right)^2\Vert z-z_0\Vert.\nonumber\\
			\intertext{Hence, for all $z\in\mathcal{Z},$}
			\vert \phi_k(z,u)-\phi_k(z_0,u)\vert \leq \left(1+\left(\gamma_{R_k}+\gamma_{D_k}\right)\right) \Vert z-z_0\Vert\phantom{iiii}\nonumber\\\times\left[\lambda \displaystyle\sum_{j=1}^n \tilde{C}_{kj}\textcolor{black}{\gamma_{D_j}N_j}+\theta_{sup}+\displaystyle\left\lvert\dfrac{\displaystyle\sum_{j=1}^n \tilde{C}_{kj}z_{3j}}{\displaystyle\sum_{j=1}^n \tilde{C}_{kj}z_{2j}}\displaystyle\right\rvert+\left(\gamma_{R_k}+\gamma_{D_k}\right)\right] \nonumber\\
			+\left(\gamma_{R_k}+\gamma_{D_k}\right)^2\Vert z-z_0\Vert.\phantom{iiiiiiiiiiiiiiiiiiiiiiiiiiiiiiiiii}\label{eq20}
		\end{align}
		The fraction term is smaller than a constant. Indeed, 
		\begin{align*}
			\displaystyle\left\lvert\dfrac{\displaystyle\sum_{j=1}^n \tilde{C}_{kj}z_{3j}}{\displaystyle\sum_{j=1}^n C_{kj}z_{2j}}\displaystyle\right\rvert\leq\dfrac{\displaystyle\sum_{j=1}^n \tilde{C}_{kj}\vert z_{3j}\vert}{\displaystyle\left\lvert\displaystyle\sum_{j=1}^n \tilde{C}_{kj}z_{2j}\displaystyle\right\rvert}\phantom{iiiiiiiiiiiiiiiiiiiiiiiiiiiiiiiiiiiiiiiiiiiiiiii}\\
			=\dfrac{\displaystyle\sum_{j=1}^n \tilde{C}_{kj}\displaystyle\left\lvert \left\{\gamma_{D_j}\lambda S_j \displaystyle\sum_{l=1}^n \tilde{C}_{jl}z_{2_l} -\left(\gamma_{R_j}+\gamma_{D_j}\right)z_{2j}\left(t\right)\right\} \displaystyle\right\rvert}{\displaystyle\sum_{j=1}^n \tilde{C}_{kj}z_{2j}}\phantom{iiiiiiiiiii}\\
			\leq \dfrac{\lambda \displaystyle\sum_{j=1}^n \tilde{C}_{kj}\gamma_{D_j}N_j\displaystyle\sum_{l=1}^n \tilde{C}_{jl}\vert z_{2_l}\vert }{\displaystyle\sum_{j=1}^n \tilde{C}_{kj}z_{2j}}+\dfrac{\Gamma \displaystyle\sum_{j=1}^n \tilde{C}_{kj}\vert z_{2j}\vert }{\displaystyle\sum_{j=1}^n \tilde{C}_{kj}z_{2j}}\phantom{iiiiiiiiiiiiiiiiiiii}\\
			\intertext{using the fact that $\gamma_{R_j}+\gamma_{D_j}\leq \Gamma$ for all $j$, where $\Gamma$ is given as in Theorem \ref{thm_pos_unsaturated_fdb}, and the fact that $0\leq S_j\leq N_j$ for all $j$ over the domain of interest. Then,}\\
			\displaystyle\left\lvert\dfrac{\displaystyle\sum_{j=1}^n \tilde{C}_{kj}z_{3j}}{\displaystyle\sum_{j=1}^n \tilde{C}_{kj}z_{2j}}\displaystyle\right\rvert\leq\dfrac{\lambda \displaystyle\sum_{j=1}^n \tilde{C}_{kj}\gamma_{D_j}N_j\displaystyle\sum_{\substack{l=1 \\ z_{2l}\neq 0}}^n \tilde{C}_{jl} z_{2_l} }{\displaystyle\sum_{j=1}^n \tilde{C}_{kj}z_{2j}}+\Gamma\phantom{iiiiiiiiiiiiii}\\
			=\lambda \displaystyle\sum_{j=1}^n \tilde{C}_{kj}\gamma_{D_j}N_j\left\{\displaystyle\sum_{\substack{l=1 \\ z_{2l}\neq 0}}^n\dfrac{\tilde{C}_{jl}}{\displaystyle\sum_{j=1}^n\tilde{C}_{kj}\dfrac{z_{2j}}{z_{2l}}}\right\}+\Gamma\phantom{iiiiiiiiiiiiii}\\
		\end{align*}
		\begin{align*}		
			=\lambda \displaystyle\sum_{j=1}^n \tilde{C}_{kj}\gamma_{D_j}N_j\Bigg\{\displaystyle\sum_{\substack{l=1 \\ z_{2l}\neq 0}}^n\dfrac{\tilde{C}_{jl}}{\tilde{C}_{kl}\left(1+\displaystyle\sum_{\substack{j=1 \\ j\neq l}}^n\dfrac{\tilde{C}_{kj}}{\tilde{C}_{kl}}\dfrac{z_{2j}}{z_{2l}}\right)}\Bigg\}+\Gamma\phantom{iiiiiiiiiiii}\\
			\leq K \phantom{iiiiiiiiiiiiiiiiiiiiiiiiiiiiiiiiiiiiiiiiiiiiiiiiiiiiiiiiiiiiiiiiiiiiiiiiiii}
		\end{align*}
		where $K=\lambda \displaystyle\sum_{j=1}^n \tilde{C}_{kj}\gamma_{D_j}N_j\displaystyle\sum_{l=1}^n\dfrac{\tilde{C}_{jl}}{\tilde{C}_{kl}}+\Gamma$. Therefore, inequality \eqref{eq20} becomes
		\begin{align*}
			\vert \phi_k(z,U)-\phi_k(z_0,U)\vert
			\leq M \Vert z-z_0\Vert 
		\end{align*}
		where $M=\left(1+\Gamma\right)\left[\lambda \displaystyle\sum_{j=1}^n \tilde{C}_{kj}\textcolor{black}{\gamma_{D_j}N_j}+\theta_{sup}+K+\Gamma\right]+\Gamma^2 >0$. \\
	 \end{proof}

The following theorem states the main result of this section.\\
	\begin{thm}\label{thm_HG}
		Consider model \eqref{state-feed} and the high-gain observer given by
\begin{equation}\label{HG}
	\dot{\hat{z}}(t)=A\hat{z}(t)+B\phi(\hat{z}(t),\textcolor{black}{u}(\hat{z}(t)))+G(y(t)-H\hat{z}(t)),
\end{equation}
with $\hat{z}(0)=\hat{z}_0$, where $G$ denotes the observer gain, defined by $G=\text{blockdiag}[G_1,...,G_n]_{3n\times n}$ where $$G_i= \begin{pmatrix}
	\dfrac{\beta_1^i}{\epsilon}&\dfrac{\beta_2^i}{\epsilon^2}&\dfrac{\beta_3^i}{\epsilon^3}
\end{pmatrix}^T$$ with the parameters $\beta_i^j$, $j=1, 2, 3$ chosen such that the roots of $s^3+\beta_1^is^2+\beta_2^is+\beta_3^i $ are in the open left-half plane, for $i=1,...,n$.
Then, 
there exists $\tilde{\epsilon}^{\star}$ such that, for every $0<\epsilon\leq \tilde{\epsilon}^{\star}$, the equilibrium set $\tilde{Z}^{\star}:=\textcolor{black}{\left\{\tilde{z}\in \mathcal{Z} : \tilde{z}_{2k}=\tilde{z}_{3k}=0, k=1,...,n\right\}}$ of system \eqref{state-feed} under observer-based feedback is asymptotically stable and pointwise attractive.
\end{thm}\vspace{0.3cm}
\begin{proof}
For the analysis part, one can replace the observer dynamics by the equivalent dynamics of the scaled estimation error: $$\eta_{ik}=\dfrac{\tilde{z}_{ik}-\hat{z}_{ik}}{\epsilon^{3-i}}, $$ for $k=1,...,n$.\\
	It follows that $\hat{z}=\tilde{z}-D(\epsilon)\eta$ where 
	\begin{align*}
		\eta&=\left(\eta_{11},\eta_{21}, \eta_{31}, \eta_{12},\cdots,\eta_{3n}\right)^T\\
		D(\epsilon)&=\text{blockdiag}[\tilde{D},...,\tilde{D}]_{3n\times 3n}\\
		\tilde{D}&=\begin{pmatrix}
		\epsilon^2 & 0 & 0\\
		0 & \epsilon & 0\\
		0 & 0 & 1
		\end{pmatrix}.	
	\end{align*}
The closed-loop system is then given by
	\begin{equation}\label{HGO_sys}
\left\{\begin{array}{rl}
\dot{\tilde{z}}\left(t\right) &= A\tilde{z}\left(t\right)+B\phi\left(\tilde{z}\left(t\right),\textcolor{black}{u}(\tilde{z}(t)-D(\epsilon)\eta\left(t\right))\right) \\ 
\epsilon\dot{\eta}\left(t\right) &=A_0\eta+\epsilon B g\left(\tilde{z}(t),\tilde{z}(t)-D(\epsilon)\eta\left(t\right)\right)
\end{array}\\
\right.
\end{equation}
where \textcolor{black}{	\begin{align*}
	A_0&=\text{blockdiag}[\tilde{A}_{0_1},...,\tilde{A}_{0_n}]_{3n\times 3n}\\
	\tilde{A}_{0_i}&=\begin{pmatrix}
		-\beta_1^i & 1 & 0\\
		-\beta_2^i & 0 & 1\\
		-\beta_3^i & 0 & 0
	\end{pmatrix}	
\end{align*}} and $g\left(\tilde{z},\tilde{z}-D(\epsilon)\eta\right)=$\\ $\phantom{and g()}\phi\left(\tilde{z},\textcolor{black}{u}(\tilde{z}-D(\epsilon)\eta)\right)-\phi\left(\tilde{z}-D(\epsilon)\eta,\textcolor{black}{u}(\tilde{z}-D(\epsilon)\eta)\right).$\vspace{0.2cm}\\
In the following, the notation $\chi=\left(\tilde{z}^T, \eta^T\right)^T$ is introduced.
The set of equilibria of model \eqref{HGO_sys} is given by $$\chi^{\star}=:\tilde{Z}^{\star}\times 0_{\mathbb{R}^{3n}}. $$
LaSalle's theorem implies the conclusion about the attractivity of the set $\chi^{\star}$. Indeed, consider the function $V:\mathcal{Z}\times \mathbb{R}^3 \to \mathbb{R}$ such that for all $\chi\in\mathcal{Z}\times \mathbb{R}^3,$ $$V(\chi)=\displaystyle\sum_{k=1}^n (N_k-\tilde{z}_{1k})+\eta^TP\eta\geq 0,$$ where $P=P^T$ is the positive definite solution of the Lyapunov equation $PA_0+A_0^TP=-I$. The time derivative of $V$ along the trajectories of model \eqref{HGO_sys} is given by 
\begin{align*}
\dot{V}(\chi)&=\displaystyle\sum_{k=1}^n -\tilde{z}_{2k}+\dot{\eta}^TP\eta+\eta^TP\dot{\eta}\\
&=\displaystyle\sum_{k=1}^n -\tilde{z}_{2k}+\dfrac{1}{\epsilon}\left[\left(A_0\eta+\epsilon Bg(\tilde{z},\tilde{z}-D(\epsilon)\eta)\right)^TP\eta\right.\\&\left.\phantom{iiiiiiii}+\eta^TP\left(A_0\eta+\epsilon Bg(\tilde{z},\tilde{z}-D(\epsilon)\eta)\right)\right]\\
&=\displaystyle\sum_{k=1}^n -\tilde{z}_{2k}+\dfrac{1}{\epsilon}\left[\eta^T\left(A_0^TP+PA_0\right)\eta\right.\\
&\left.\phantom{iiiiiiii}+2\epsilon\langle P\eta, Bg(\tilde{z},\tilde{z}-D(\epsilon\eta))\rangle\right]\\
&\leq \displaystyle\sum_{k=1}^n -\tilde{z}_{2k}+\dfrac{1}{\epsilon}\left[-\Vert\eta\Vert^2
+2\epsilon\Vert P\eta\Vert\hspace{0.1cm} \Vert Bg(\tilde{z},\tilde{z}-D(\epsilon\eta))\Vert\right],\\
\intertext{thanks to Cauchy-Schwarz inequality. \textcolor{black}{Moreover, in view of Remark \ref{rem_lip}}, Lemma \ref{lemma_ass_1} can be used to get the following inequalities:}
\dot{V}(\chi)&\leq \displaystyle\sum_{k=1}^n -\tilde{z}_{2k}+\dfrac{1}{\epsilon}\left[-\Vert\eta\Vert^2+2\epsilon\Vert P\eta\Vert\hspace{0.1cm} M \Vert D(\epsilon)\eta \Vert \right]\\
&\leq \displaystyle\sum_{k=1}^n -\tilde{z}_{2k}+\left[\dfrac{-1}{\epsilon}+2\Vert P\Vert\hspace{0.1cm}  M \Vert D(\epsilon) \Vert\right] \Vert\eta\Vert^2\\
&\leq 0 
\end{align*}
for all $0<\epsilon<\tilde{\epsilon}^{\star}$ where $\tilde{\epsilon}^{\star}$ is chosen such that  $\tilde{\epsilon}^{\star}D(\tilde{\epsilon}^{\star})= \dfrac{1}{2\Vert P  \Vert}.$\\ Moreover, one can notice that $\dot{V}(\chi)=0\Leftrightarrow \tilde{z}_{2k} (\Rightarrow \tilde{z}_{3k})=0$ and $\eta=0$. Hence, by LaSalle's theorem, any solution starting in $Z\times \mathbb{R}^3$ converges to the set $\chi^{\star}$. Furthermore, in view of the dynamics of $\tilde{z}_{1k}$ and the nonnegativity of $\tilde{z}_{2k}$ in $Z$, it follows that $\tilde{z}_{1k}$ is increasing and bounded, hence it converges to a point $\tilde{z}_{1k}^{\star}$. Therefore, $\chi$ converges to a point in $\chi^{\star}$. The definition of $\chi^{\star}$ concludes the proof.

\end{proof}
	
	\section{Numerical simulations}\label{Sec_num_sim}
	In this section, numerical simulations, illustrating Propositions \ref{prop_asympt_stab} and \ref{thm_stab_sat_vacc} and Theorems \ref{thm_pos_unsaturated_fdb} and \ref{thm_HG}, are presented. Those simulations consist of academic examples where the parameters were chosen thoroughly to be close to the reality of some given epidemic.  The model is divided in $6$ classes of age, $[0-29), [29-39), [39-49), [49-59), [59-90), 90+$, since  real data are taken for the social contact matrix from Socrates, an online tools from \cite{willem_2020}, \cite{social_contact_data}. For these simulations, data about the first covid wave in France are chosen. Even though the authors are aware that an SIRD model with no loss of immunity is not the best choice to model COVID-19, parameters coming from a real epidemic give a good illustrative example to the theory developped here. The choice of the recovery rate and the death rate is motivated by the work of \cite{guan_2020}, where, for the region of Provence-Alpes-Côte d'Azur, those rates equal to $0.1722$ and $0.0242$ respectively. In the following simulations, the parameters are taken with the same order of magnitude but emphasize the fact that older people are more difficult to heal and more likely to die from the disease. This leads to the following choice, $\gamma_R=\begin{pmatrix}
		0.3& 0.3& 0.3& 0.1& 0.1& 0.1
	\end{pmatrix}^T$, $\gamma_D=\begin{pmatrix}0.001& 0.01& 0.01& 0.04& 0.05& 0.15\end{pmatrix}^T$. The disease transmission probability is given by $\lambda=0.5$. This means that the probability to catch the disease is equal to $50\%$. Furthermore, the immunization is assumed to work perfectly, so $p_k=1$ for $k=1,...,n$. Moreover, real data, coming from \cite{insee} are taken for the distribution of the number of individuals by class of age, $N_k$, $k=1,...,6$ where we focus on the Vaucluse department ($84$). Concerning the initial conditions, it is assumed there are only a few infected individuals in some classes of the population at the beginning of the disease: $I_0=\begin{pmatrix}
		0& 0& 20& 30& 0& 0
	\end{pmatrix}$ and that there is no recovered people and deceased people yet, since it is the start of the epidemic. Hence, the first susceptible individuals are computed using the following relation: $S_k(0)=N_k-I_k(0)$. Finally, the simulation is stopped when some convergence criterion is satisfied. In this case, it corresponds to the time when there remains less than one infected individuals in the whole population.
	\subsection{Open-loop case}
\begin{figure}[t]
		\centering
		\includegraphics[height=5cm]{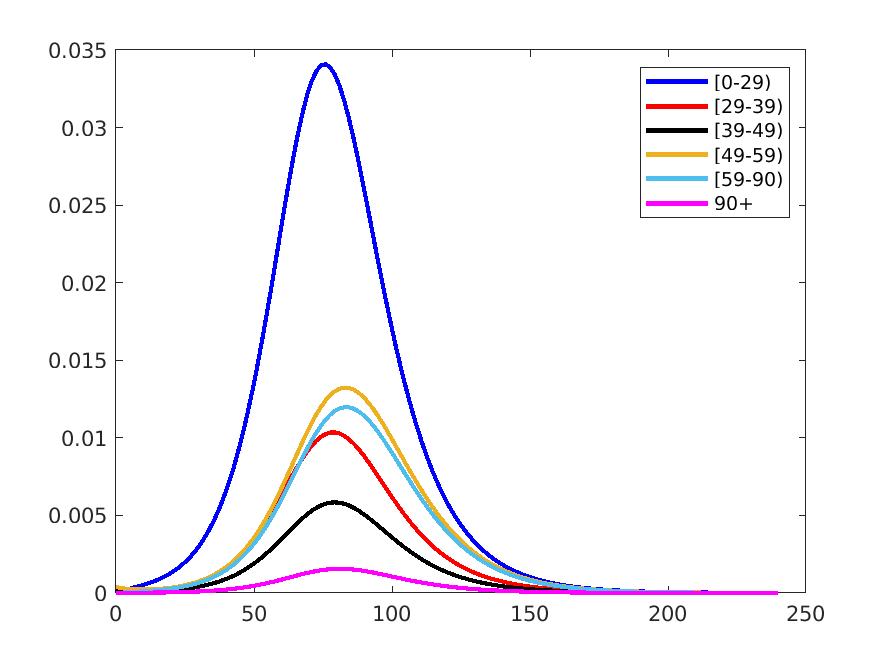}
		\caption{Dynamics of the proportion of infected individuals without control}  
		\label{fig1}   
	\end{figure}

\begin{figure}
		\centering
		\includegraphics[height=5cm]{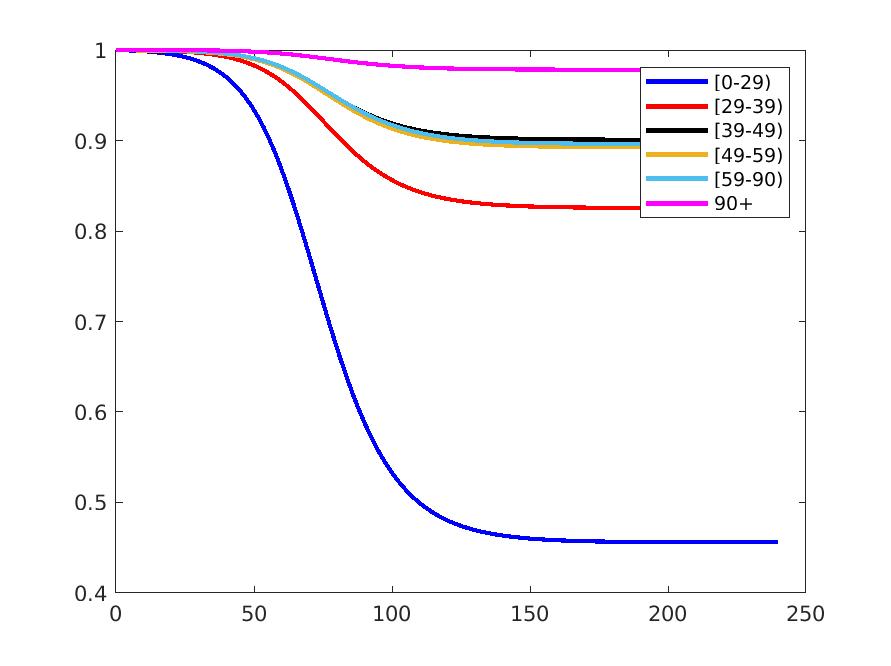}
		\caption{Dynamics of the proportion of susceptible individuals without control}
		\label{fig2}   
		
	\end{figure}

	This part is dedicated to the numerical simulation in the open-loop case. Figures \ref{fig1} and \ref{fig2} show the dynamics of the infected and susceptible individuals, respectively, without control. In Figure \ref{fig1}, the convergence of the infected individuals to zero is shown. It happens after $240$ days. This means that, according to the chosen SIRD model, this hypothetical epidemic will disappear after $240$ days as there are no longer any infected individuals. Moreover, the higher number of infected individuals appears in the first class of age, despite the fact that there are no infected individuals in this class of age at the beginning of the disease. This is probably due to the high contact rate of this class of age and the higher number of individuals in this class. To avoid hospital burden,  the objective is to minimize this number, accounting for all age groups. This justifies the introduction of a feedback law, which will be further discussed in the following subsections. Moreover, We can observe on Figure \ref{fig2}, that the susceptible individuals also tend to an equilibrium, which differs for each class of age as predicted by Proposition \ref{prop_asympt_stab}. 	
	
	\subsection{Unsaturated state feedback}
\begin{figure}[t]
		\centering
		\includegraphics[height=5cm]{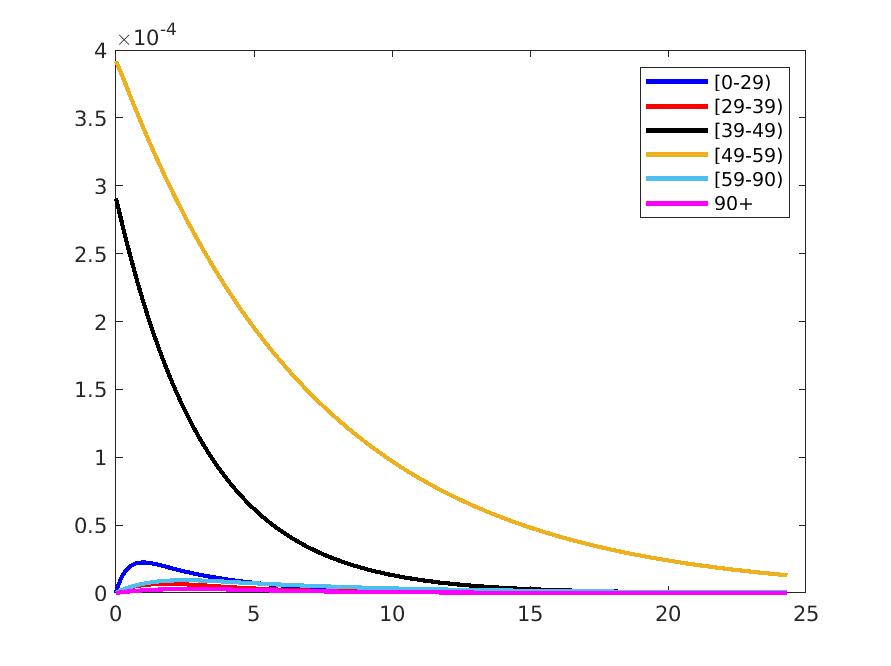}
		\caption{Dynamics of the proportion of infected individuals with unconstrained state feedback}  
		\label{fig3} 
	\end{figure}

	In view of Figure \ref{fig1}, one can wish to obtain disease eradication faster and with less infected individuals in the population. This can be done with the state feedback law \eqref{control} that implies exponential convergence to zero of the infected individuals. This law was implemented satisfying parameter conditions mentioned in Theorem \ref{thm_pos_unsaturated_fdb}, with $\alpha_1^k=(\gamma_{R_k}+\gamma_{D_k})\left(max(\gamma_{D_k}+\gamma_{R_k})+\displaystyle\sum_{j=1}^n M_{kj}\right)+0.1$ and $
	\alpha_2^k=\gamma_{R_k}+\gamma_{D_k}+max(\gamma_{D_k}+\gamma_{R_k})+\displaystyle\sum_{j=1}^n M_{kj}.$ \\Figure \ref{fig3} shows that the convergence of individuals to zero is much faster than without control. This occurs in $25$ days compared to $240$ days in absence of control. Moreover, the peak of infected individuals is much lower. By focusing on the first class of age, the proportion of infected individuals is at most $0.266\times 10^{-4}$ much smaller than $0.0341$ in Figure~\ref{fig1}. Remark that this phenomenon of reduction is also seen for the dead individuals, as it can be expected. Indeed, in simulation, whose graphs are not presented here, a proportion of $0.018$ deceased individuals is obtained in the open-loop case, for the first class of age, but only $9.85 \times 10^{-8}$ when the control law is applied. Finally, the dynamics of the control law, not represented here, tends quickly to infinity as expected and will not be applicable in practice. Thus, the state feedback law \eqref{control} enables fast disease eradication but is not implementable.
	\subsection{Amplitude constrained control}
\begin{figure}
		\centering
		\includegraphics[height=5cm]{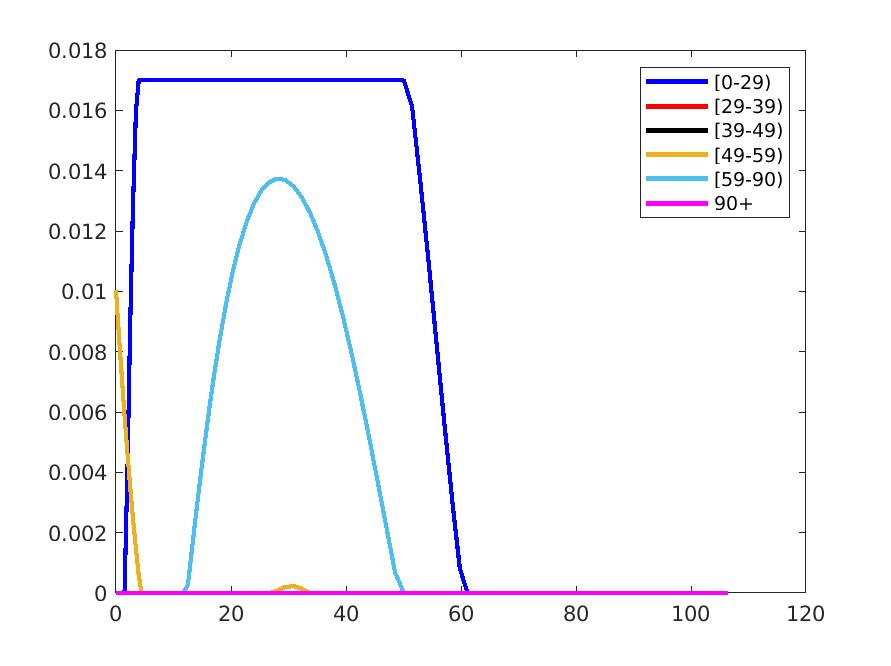}
		\caption{Dynamics of the constrained control law}  
		\label{fig6} 
	\end{figure}

\begin{figure}
	\centering
	\includegraphics[height=5cm]{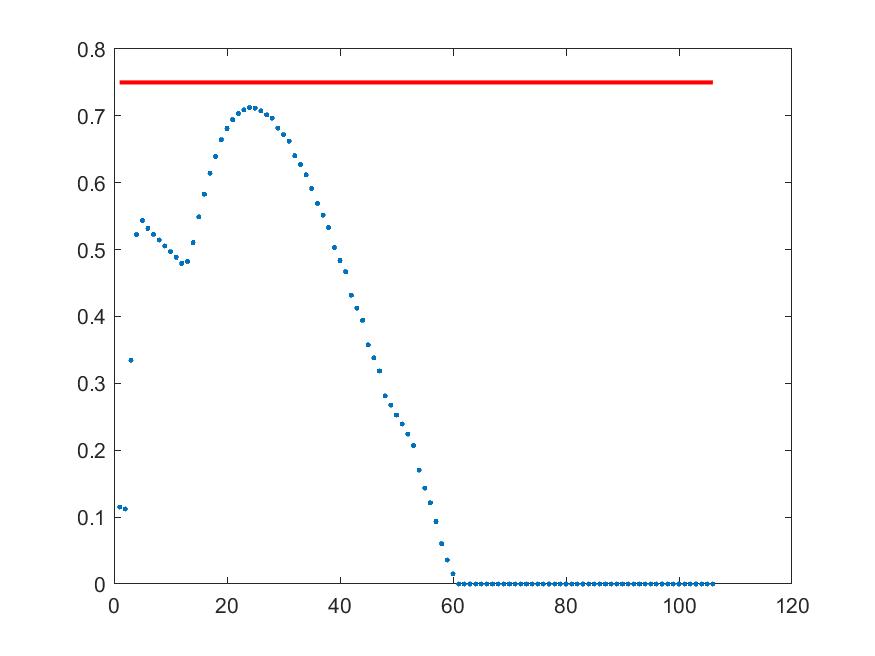}
	\caption{Percentage of vaccinated people per day}  
	\label{fig9} 
\end{figure}

\begin{figure}
	\centering
	\includegraphics[height=5cm]{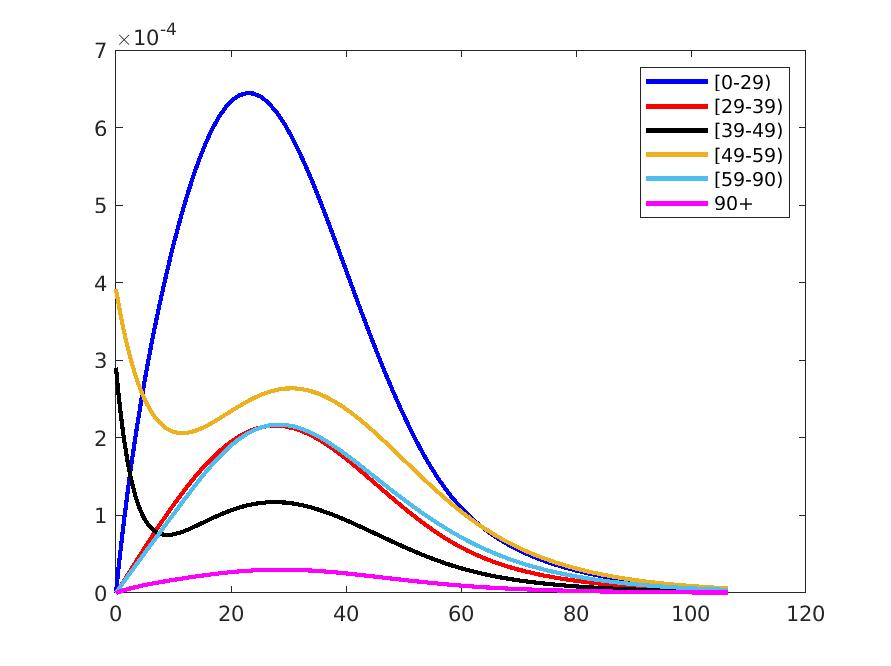}
	\caption{Dynamics of the proportion of infected individuals with constrained state feedback}  
	\label{fig7} 
\end{figure}

\begin{figure}
	\centering
	\includegraphics[height=5cm]{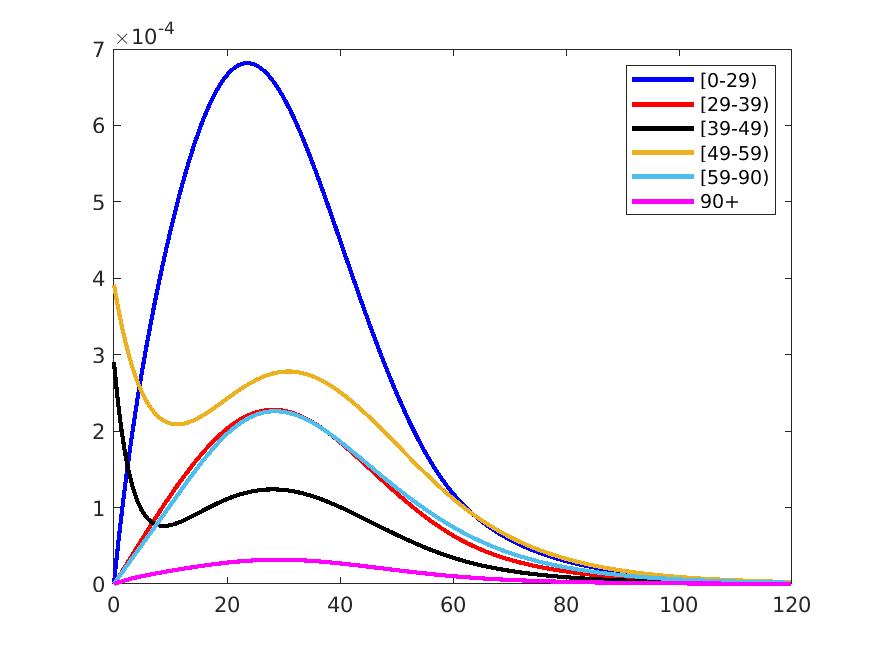}
	\caption{Dynamics of the proportion of infected individuals with observer-based output feedback, for $\epsilon=0.01$}  
	\label{fig8} 
\end{figure}

	To ensure a feasible control law, the state feedback is saturated and is given by \eqref{control_sat}. The upper bound is fixed to $0.017$ to ensure a physical meaning, as shown in Figure \ref{fig6}. 
	This value has been determined through trial and error to achieve a vaccination rate of $0.75\%$ of the population per day, serving as a threshold. This is depicted in  Figure \ref{fig9} illustrating that under the assumption that this maximum is satisfied, the proposed control law can be practically implemented.
 Moreover, the switch in the control law $\theta_{sat_k}(t)$ is chosen to satisfy Lemma \ref{lemma_invariant} where $\tilde{I}_k$ is arbitrarily set to $20$ for all $k$ and the equality is taken to choose $\tilde{S}_k$. Under this feedback, the convergence of the infected individuals to $0$ is no longer exponential as it can be observed in Figure \ref{fig7}. We have disease eradication with a smaller peak of infected individuals than in the open-loop case, a proportion of $0.0015$ instead of $0.0759$. The same holds for the dead individuals (not represented here), a proportion of $0.0015$ is observed with the feedback instead of $0.0891$ in open-loop. Moreover, it is interesting to observe on Figure \ref{fig6} that the law is highly age-dependent. It recommends to focus effort on the class of age $[0-29)$, $[59-90)$ and a bit on $[49-59)$, and it is also time dependent. At first, it suggests to begin by the age between $49$ and $59$. Moreover, this example illustrates the fact that the control can be switched on and off a finite number of time. For instance, if we focus on the class of age $[49-59)$, control is needed, then the proportion of infected individuals goes below the threshold, so the control is stopped for this class of age. Then, the proportion goes above the threshold and the control needs to be reactivate for this class of age between days $28$ and $33$. 
Moreover, since the input $\theta_{sat_k}$ is a state feedback, it is mandatory to design an observer-based output feedback, which could be used in practice. Theorem \ref{thm_HG} shows that this law can be implemented thanks to an observer-based output feedback and recovers the performance of $\theta_{sat_k}$, in the sense of Proposition \ref{prop_asympt_stab}. This is illustrated in Figure \ref{fig8}, where the initial conditions for the estimated states are given by $\hat{I}_{0}=\begin{pmatrix}
			100& 100& 100& 100& 100& 100
		\end{pmatrix}$ and $\hat{S}_{k0}=N_k-\hat{I}_{k0}$. One can observe that the estimated closed-loop trajectory of the infected individuals converges to the disease-free equilibrium. In this case, the peak of infected individuals is a little bit greater than in the case with state-feedback control law, a maximal total proportion of $0.00178$ infected individuals instead of $0.0015$. However, it is still lower than in the open-loop case, which has a maximal total proportion of $0.0759$ infected individuals.  
Hence, relying solely on data concerning deceased individuals, an age-dependent vaccination policy is implemented. This policy looks like the state feedback law introduced in Figure \ref{fig6}, although it is not explicitly depicted. As for Figure \ref{fig6}, its interpretation offers a comprehensive insight into which age groups to prioritize for vaccination, aiming to minimize the peak of infected individuals within the population.
	


	\section{Conclusion and perspectives}
	The dynamical analysis of an age-structured SIRD model was performed and results about pointwise asymptotic stability were obtained. In order to improve disease eradication in the sense of decreasing the peak of infected individuals in the population, a state feedback law was designed and exponential convergence properties were obtained. However, since this law was not implementable in practice, a saturated sate-feedback law was designed and shown to  decrease the peak of infected individuals. However, state feedback cannot be used in practice, therefore an observer-based output feedback law was obtained.\\ \\ Some complementary research will still be relevant. For instance, since there exists an infinity of state feedback gains, one can try to find those parameters under some given conditions. For instance, find the parameters that minimize the number of dead people instead of just focusing on the infected individuals, as it is done in this paper. It would also be interesting to find a feedback that satisfies some constraints on the state given a priori. One can cite the problem of finding a feedback that prevents the maximum number of infected individuals to be larger than a factor of the hospital beds capacity. In addition, while the strategy under consideration is assumed to confer long-term immunity, it is important to acknowledge that boosters may be necessary in certain cases to sustain immunity over time. This scenario could be explored by permitting recovered individuals to transition back to the susceptible compartment, resulting in a modified model known as SIRDS. Moreover, it could be very interesting to perform the numerical simulation on real data that fits the SIRD model developed in this work. Therefore, a parameter identification approach needs to be implemented in order to find the appropriate model parameters to represent a relevant epidemic. Considering the necessity of using real data, it is crucial to note that the measurements are sampled. Therefore, it would be interesting to extend Section \ref{Sect_output_fdb} to the case of sampled-data high-gain observers. The article \cite{Dabroom_2001} could offer valuable insights to address this question. Moreover, modifications of model \eqref{SIRD} could be studied. Those modifications can occur by adding other compartments to the model such as the compartment $E$ of exposed individuals for instance. It could also be interesting to incorporate the aging effect in the dynamics to focus on long-term diseases. Finally, it could be interesting to compare the results obtained in this paper to some obtained with different control strategy.

	

	\appendix[Lipschitz property of the saturated feedback]
	This section is dedicated to the proof of intermediate tools useful for the proof of the Lipschitz property of the constrained feedback \eqref{control_sat}. The areas mentioned in the proof refer to Figure \ref{fig_lipschitz}.\vspace{0.3cm}
	\begin{lemma}
		The function $q_k:Dom(q_k)=\left[\tilde{S}_k,N_k\right]\times\left[0,N_k\right]\times \cdots \times \left[0,N_k\right] \times \left[\tilde{I}_k,N_k\right]\times \left[0,N_k\right]\times \cdots \times \left[0,N_k\right]\to \mathbb{R}$, defined by \eqref{fct_gk} is Lipschitz. 
	\end{lemma}\vspace{0.3cm}
	\begin{proof}
		Let $z_k=(x_k,y_1,\cdots,y_k,\cdots,y_n)^T$ and $z_k'=(x_k',y'_1,\cdots,y'_k,\cdots,y'_n)^T \in Dom(q_k)$. Several situations can be identified. 
		\begin{itemize}
			\item Assume that $z_k$ and $z_k'$ are such that $x_k, x_k'\geq \tilde{\tilde{S}}_k$ and $y_k, y_k'\geq \tilde{\tilde{I}}_k$ ($z_k$ and $ z_k'\in B_1$). Therefore, $$\vert q_k(z_k)-q_k(z_k') \vert =\vert 1-1 \vert =0 \leq C_1 \Vert z_k-z_k'\Vert$$ for all $C_1>0$.
			\item Assume that $z_k$ and $z_k'$ are such that $y_k, y_k'< \tilde{\tilde{I}}_k$, $y_k\leq \dfrac{ \tilde{\tilde{I}}_k-\tilde{I}_k}{\tilde{\tilde{S}}_k-\tilde{S}_k}\left(x_k-\tilde{S}_k\right)+\tilde{I}_k$ and $y_k'\leq \dfrac{ \tilde{\tilde{I}}_k-\tilde{I}_k}{\tilde{\tilde{S}}_k-\tilde{S}_k}\left(x_k'-\tilde{S}_k\right)+\tilde{I}_k$ ($z_k$ and $z_k'\in B_2$). Therefore, \begin{align*}
				\left| q_k(z_k)-q_k(z_k') \right| &=\dfrac{4}{\pi}\left| \arctan\left(\dfrac{y_k-\tilde{I}_k}{\tilde{\tilde{I}}_k-\tilde{I}_k}\right) \right.\\&\hspace{0.5cm}\left.- \arctan\left(\dfrac{y_k'-\tilde{I}_k}{\tilde{\tilde{I}}_k-\tilde{I}_k}\right)\right|\\&\leq \dfrac{4}{\pi}\left| \left(\dfrac{y_k-\tilde{I}_k}{\tilde{\tilde{I}}_k-\tilde{I}_k}\right) - \left(\dfrac{y_k'-\tilde{I}_k}{\tilde{\tilde{I}}_k-\tilde{I}_k}\right)\right|\\&=\dfrac{4}{\pi\left(\tilde{\tilde{I}}_k-\tilde{I}_k\right)} \vert y_k-y_k'\vert \\ &\leq C_2 \Vert z_k-z_k'\Vert,
			\end{align*}
			using the fact that $\arctan$ is a Lipschitz function with Lipschitz constant $1$, and taking $C_2$ equals $\dfrac{4}{\pi\left(\tilde{\tilde{I}}_k-\tilde{I}_k\right)}$.
			\item Assume that $z_k$ and $z_k'$ are such that $x_k, x_k'< \tilde{\tilde{S}}_k$, $y_k> \dfrac{ \tilde{\tilde{I}}_k-\tilde{I}_k}{\tilde{\tilde{S}}_k-\tilde{S}_k}\left(x_k-\tilde{S}_k\right)+\tilde{I}_k$ and $y_k'> \dfrac{ \tilde{\tilde{I}}_k-\tilde{I}_k}{\tilde{\tilde{S}}_k-\tilde{S}_k}\left(x_k'-\tilde{S}_k\right)+\tilde{I}_k$ ($z_k$ and $z_k'\in B_3$). By a similar reasoning as the previous case, it follows that $$\left| q_k(z_k)-q_k(z_k') \right|\leq C_3 \Vert z_k-z_k'\Vert,$$ with $C_3=\dfrac{4}{\pi\left(\tilde{\tilde{S}}_k-\tilde{S}_k\right)}$.
			\item  Assume that $z_k$ and $z_k'$ are such that $x_k\geq \tilde{\tilde{S}}_k$, $y_k\geq  \tilde{\tilde{I}}_k$, $y_k'< \tilde{\tilde{I}}_k$ and $y_k'\leq \dfrac{ \tilde{\tilde{I}}_k-\tilde{I}_k}{\tilde{\tilde{S}}_k-\tilde{S}_k}\left(x_k'-\tilde{S}_k\right)+\tilde{I}_k$ ($z_k \in B_1$ and $z_k' \in B_2$). Therefore, \begin{align*}
				\left| q_k(z_k)-q_k(z_k') \right| &=\left|1-\dfrac{4}{\pi} \arctan\left(\dfrac{y_k'-\tilde{I}_k}{\tilde{\tilde{I}}_k-\tilde{I}_k}\right)\right|\\&=\left|\dfrac{4}{\pi} \arctan\left(1\right)-\dfrac{4}{\pi} \arctan\left(\dfrac{y_k'-\tilde{I}_k}{\tilde{\tilde{I}}_k-\tilde{I}_k}\right)\right|\\&\leq \dfrac{4}{\pi}\left|1 - \left(\dfrac{y_k'-\tilde{I}_k}{\tilde{\tilde{I}}_k-\tilde{I}_k}\right)\right|\\&=\dfrac{4}{\pi\left(\tilde{\tilde{I}}_k-\tilde{I}_k\right)} \left( \tilde{\tilde{I}}_k-y_k'\right) \text{ since } y_k'< \tilde{\tilde{I}}_k,\\&=C_2\left(y_k-y_k'+\tilde{\tilde{I}}_k-y_k\right) \\ &\leq C_2\left(y_k-y_k'\right)\text{ since } y_k\geq \tilde{\tilde{I}}_k,\\ &\leq C_2 \Vert z_k-z_k'\Vert,
			\end{align*}
			\item Assume that $z_k$ and $z_k'$ are such that $x_k\geq \tilde{\tilde{S}}_k$, $y_k\geq  \tilde{\tilde{I}}_k$, $x_k'< \tilde{\tilde{S}}_k$ and $y_k'> \dfrac{ \tilde{\tilde{I}}_k-\tilde{I}_k}{\tilde{\tilde{S}}_k-\tilde{S}_k}\left(x_k'-\tilde{S}_k\right)+\tilde{I}_k$ ($z_k \in B_1$ and $z_k' \in B_3$). A similar reasoning as the previous one gives $$\left| q_k(z_k)-q_k(z_k') \right|\leq C_3 \Vert z_k-z_k'\Vert.$$
			\item Assume that $z_k$ and $z_k'$ are such that $y_k< \tilde{\tilde{I}}_k$ and $y_k\leq \dfrac{ \tilde{\tilde{I}}_k-\tilde{I}_k}{\tilde{\tilde{S}}_k-\tilde{S}_k}\left(x_k-\tilde{S}_k\right)+\tilde{I}_k$, $x_k'< \tilde{\tilde{S}}_k$ and $y_k'> \dfrac{ \tilde{\tilde{I}}_k-\tilde{I}_k}{\tilde{\tilde{S}}_k-\tilde{S}_k}\left(x_k'-\tilde{S}_k\right)+\tilde{I}_k$ ($z_k \in B_2$ and $z_k' \in B_3$).
			\begin{align*}
				\left| q_k(z_k)-q_k(z_k') \right| &=\dfrac{4}{\pi}\left|\arctan\left(\dfrac{y_k-\tilde{I}_k}{\tilde{\tilde{I}}_k-\tilde{I}_k}\right)\right.\\&\hspace{0.5cm}\left.- \arctan\left(\dfrac{x_k'-\tilde{S}_k}{\tilde{\tilde{S}}_k-\tilde{S}_k}\right)\right|\\&\leq \dfrac{4}{\pi}\left|\left(\dfrac{y_k-\tilde{I}_k}{\tilde{\tilde{I}}_k-\tilde{I}_k}\right) - \left(\dfrac{x_k'-\tilde{S}_k}{\tilde{\tilde{S}}_k-\tilde{S}_k}\right)\right|
			\end{align*}
			Introduce the new variables $\bar{I}_k=\tilde{\tilde{I}}_k-\tilde{I}_k$ and $\bar{S}_k=\tilde{\tilde{S}}_k-\tilde{S}_k$ and
			consider two cases. If $\dfrac{y_k-\tilde{\tilde{I}}_k}{\bar{I}_k}\leq \dfrac{x_k'-\tilde{\tilde{S}}_k}{\bar{S}_k},$  it follows that, 
			\begin{align*}
				\left| q_k(z_k)-q_k(z_k') \right| &\leq \dfrac{4}{\pi} \left(\dfrac{x_k'-\tilde{\tilde{S}}_k}{\bar{S}_k}-\dfrac{y_k-\tilde{\tilde{I}}_k}{\bar{I}_k}\right)\\ & \leq \dfrac{4}{\pi}\left(\dfrac{y_k'-\tilde{\tilde{I}}_k}{\bar{I}_k}-\dfrac{y_k-\tilde{\tilde{I}}_k}{\bar{I}_k}\right)
				\\&=C_2\left(y_k'-y_k\right) \\  &\leq C_2 \Vert z_k-z_k'\Vert,
			\end{align*}
			where the fact that $y_k'> \dfrac{ \bar{I}_k}{\bar{S}_k}\left(x_k'-\tilde{\tilde{S}}_k-\bar{S}_k\right)+\tilde{\tilde{I}}_k-\bar{I}_k\Leftrightarrow \dfrac{y_k'-\tilde{\tilde{I}}_k}{\bar{I}_k}>\dfrac{x_k'-\tilde{\tilde{S}}_k}{\bar{S}_k}$ was used. \\Moreover, if  $\dfrac{y_k-\tilde{\tilde{I}}_k}{\bar{I}_k}> \dfrac{x_k'-\tilde{\tilde{S}}_k}{\bar{S}_k},$  it follows that, 
			\begin{align*}
				\left| q_k(z_k)-q_k(z_k') \right| &\leq \dfrac{4}{\pi} \left(\dfrac{y_k-\tilde{\tilde{I}}_k}{\bar{I}_k}-\dfrac{x_k'-\tilde{\tilde{S}}_k}{\bar{S}_k}\right)\\ & \leq \dfrac{4}{\pi}\left(\dfrac{x_k-\tilde{\tilde{S}}_k}{\bar{S}_k}-\dfrac{x_k'-\tilde{\tilde{S}}_k}{\bar{S}_k}\right)
				\\&=C_3\left(x_k-x_k'\right) \\  &\leq C_3 \Vert z_k-z_k'\Vert,
			\end{align*}
			where the fact that $y_k\leq \dfrac{ \bar{I}_k}{\bar{S}_k}\left(x_k-\tilde{\tilde{S}}_k-\bar{S}_k\right)+\tilde{\tilde{I}}_k-\bar{I}_k\Leftrightarrow \dfrac{y_k-\tilde{\tilde{I}}_k}{\bar{I}_k}\leq\dfrac{x_k-\tilde{\tilde{S}}_k}{\bar{S}_k}$ was used.
		\end{itemize}
		Finally, $$\left| q_k(z_k)-q_k(z_k') \right| \leq C \Vert z_k-z_k'\Vert,$$
		with $C=\max\left\{C_1,C_2,C_3\right\}>0$.
	\end{proof}\vspace{0.3cm}
	\begin{lemma}
		The function $\bar{\theta}_k :Dom(\bar{\theta}_k)=\left[\tilde{S}_k,N_k\right]\times\left[0,N_k\right]\times \cdots \times \left[0,N_k\right] \times \left[\tilde{I}_k,N_k\right]\times \left[0,N_k\right]\times \cdots \times \left[0,N_k\right]\to \mathbb{R}$, defined by \eqref{fct_bar_theta_k} is Lipschitz. 
	\end{lemma}\vspace{0.3cm}
	\begin{proof}
		Let $z_k$ and $z_k'\in Dom(\bar{\theta}_k)$. As for the previous proof, several situations can be identified. 	
		\begin{itemize}
			\item Assume that $z_k$ and $z_k'$ are such that $\theta_k(z_k), \theta_k(z_k')\leq \theta_{sup}$, it follows that $$\left|\bar{\theta}_k(z_k)-\bar{\theta}_k(z_k') \right|=\left|\theta_k(z_k)-\theta_k(z_k') \right|\leq K_1\Vert z_k-z_k' \Vert $$
			since $\theta_k$, defined by \eqref{control}, is Lipschitz on $Dom(\bar{\theta}_k)$.
			\item Assume that $z_k$ and $z_k'$ are such that $\theta_k(z_k), \theta_k(z_k')< \theta_{sup}$, then $$\left|\bar{\theta}_k(z_k)-\bar{\theta}_k(z_k') \right|=\left|\theta_{sup}-\theta_{sup} \right|=0\leq K_2\Vert z_k-z_k' \Vert, $$ for all $K_2>0$.
			\item Assume that $z_k$ and $z_k'$ are such that $\theta_k(z_k)\leq \theta_{sup}$ and $\theta_k(z_k')>\theta_{sup}$. Hence,
			\begin{align*}
				\left|\bar{\theta}_k(z_k)-\bar{\theta}_k(z_k') \right|&=\left|\theta_k\left(z_k\right)-\theta_{sup} \right|\\
				&=\theta_{sup}-\theta_k\left(z_k\right) \text{ since } \theta_k\left(z_k\right)\leq \theta_{sup},\\
				&\leq \theta_k\left(z_k'\right)-\theta_k\left(z_k\right) \text{ since } \theta_k\left(z_k'\right)< \theta_{sup}\\
				&\leq K_1 \Vert z_k-z_k' \Vert
			\end{align*}
			where $K_1$ is the Lipschitz constant of $\theta_k$.
		\end{itemize}
		Thus, for all $z_k$ and $z_k'\in Dom(\bar{\theta}_k)$, $$\left|\bar{\theta}_k(z_k)-\bar{\theta}_k(z_k') \right|\leq K \Vert z_k-z_k' \Vert,$$
		with $K=\max\left\{K_1,K_2\right\}$.
	\end{proof}\\ \\
	\begin{prop}
		The function $q_k\bar{\theta}_k :Dom(q_k\bar{\theta}_k)=\left[\tilde{S}_k,N_k\right]\times\left[0,N_k\right]\times \cdots \times \left[0,N_k\right] \times \left[\tilde{I}_k,N_k\right]\times \left[0,N_k\right]\times \cdots \times \left[0,N_k\right]\to \mathbb{R}$ is Lipschitz, with Lipschitz constant $L>0$. 
	\end{prop}
	\begin{proof}
		Let $z_k$ and $z_k' \in Dom(q_k\bar{\theta}_k)$,\begin{align*}
			\left|\left(q_k\bar{\theta}_k\right)\left(z_k\right)-\left(q_k\bar{\theta}_k\right)\left(z_k'\right)\right|&\leq \left|q_k\left(z_k\right)\left(\bar{\theta}_k\left(z_k\right)-\bar{\theta}_k\left(z_k'\right)\right) \right| \\ &\hspace{0.5cm}+ \left|\bar{\theta}_k\left(z_k'\right)\left(q_k\left(z_k\right)-q_k\left(z_k'\right)\right) \right|\\ &\leq \left|\bar{\theta}_k\left(z_k\right)-\bar{\theta}_k\left(z_k'\right)\ \right|\\ &\hspace{0.5cm}+\theta_{sup} \left|q_k\left(z_k\right)-q_k\left(z_k'\right) \right|\\
			&\leq K \Vert z_k-z_k' \Vert + \theta_{sup} C \Vert z_k-z_k' \Vert\\
			&\leq L \Vert z_k-z_k' \Vert
		\end{align*}	
		using the two previous lemma, the fact that $0\leq q_k\leq 1$, $0\leq \bar{\theta}_k\leq \theta_{sup}$ and defining $L=\max\left\{K,\theta_{sup}C\right\}$.
	\end{proof}

	\ifCLASSOPTIONcaptionsoff
	\newpage
	\fi

	
	
	\bibliographystyle{IEEEtranS}
	\nocite{*}
	\bibliography{IEEEabrv,biblio}
	
	\begin{IEEEbiography}[{\includegraphics[width=1in,height=1.25in,clip,keepaspectratio]{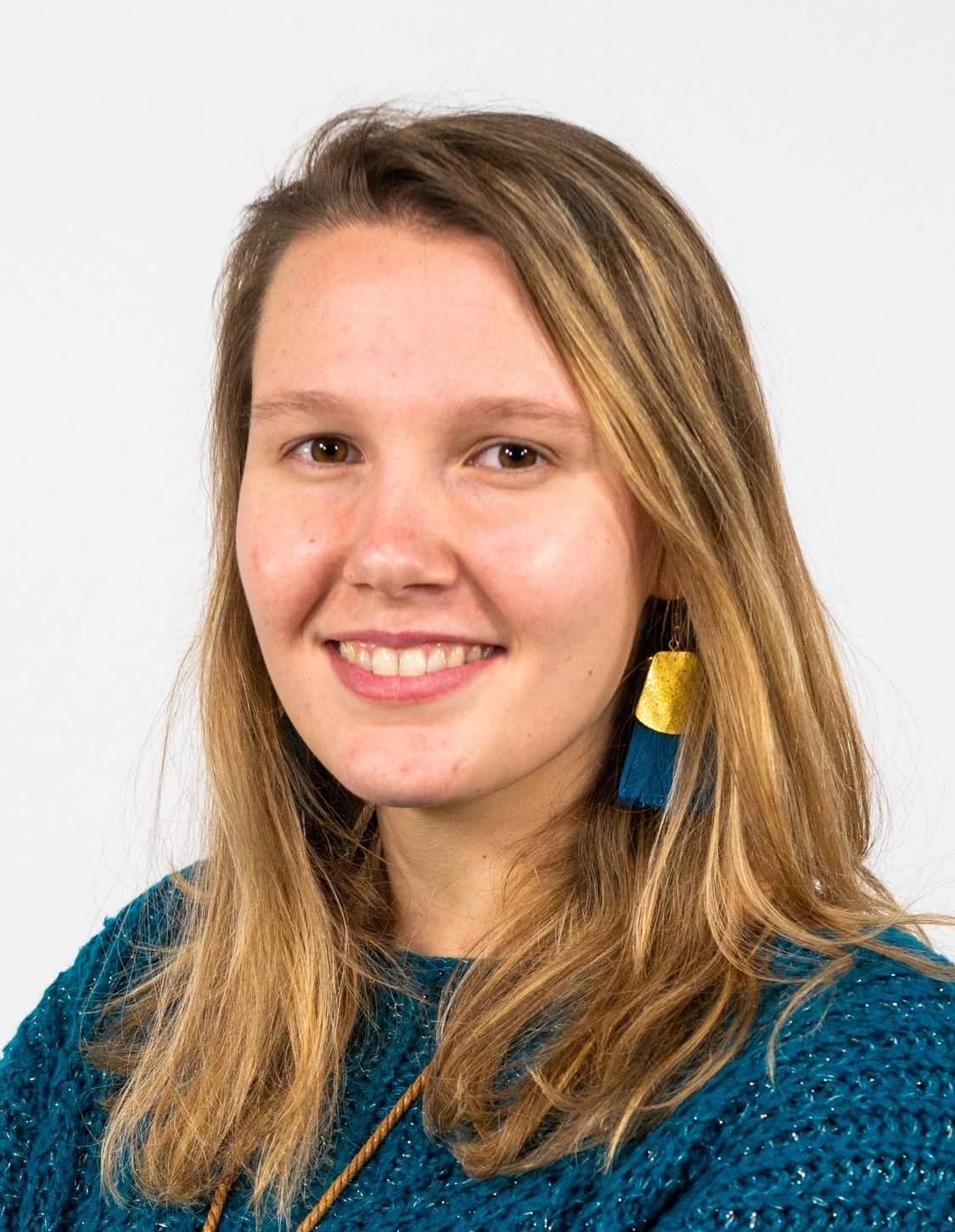}}]{Candy Sonveaux} was born in Belgium in 1995. She studied applied mathematics at the university of Namur (Belgium) and obtained her Master thesis in 2018. Currently she is a teaching assistant and a PhD Student at the department of mathematics of the university of Namur .  Her main research interest focuses on applications of system and control theory to the field of epidemic models. In particular, she is interested in dynamical analysis and control of finite-dimensional and infinite-dimensional nonlinear descriptions of such models.
\end{IEEEbiography}
	
		\begin{IEEEbiography}[{\includegraphics[width=1in,height=1.25in,clip,keepaspectratio]{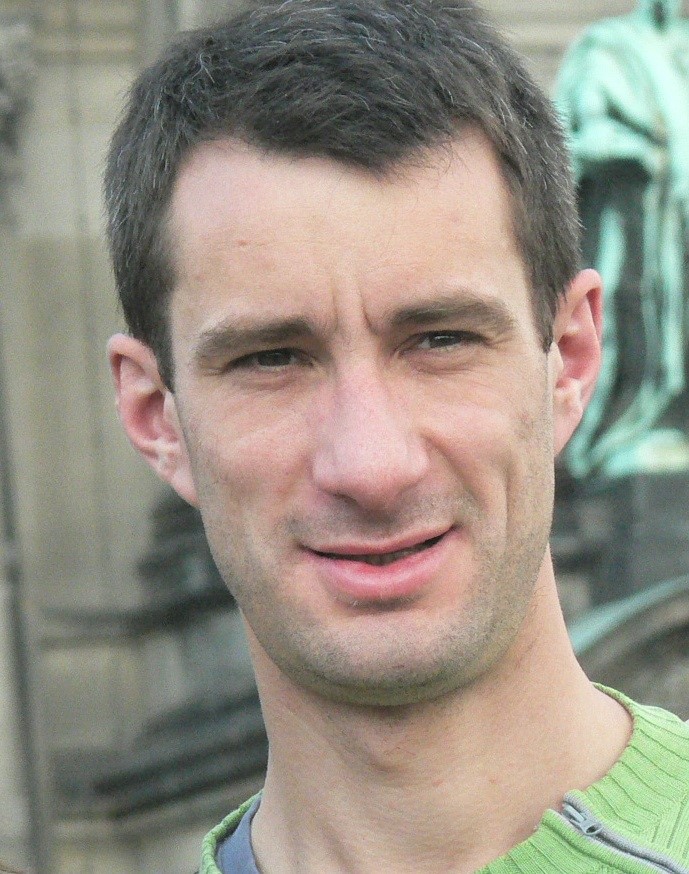}}]{Christophe Prieur} graduated in Mathematics from the Ecole Normale Supérieure de Cachan, France in 2000. He received the Ph.D degree in 2001 in Applied Mathematics from the Université Paris-Sud, France. He is currently a senior researcher of the CNRS. He is currently an associate editor of the AIMS Evolution Equations and Control Theory, the SIAM Journal of Control and Optimization and the Mathematics of Control, Signals, and Systems. He is a senior editor of the IEEE Control Systems Letters, and an editor of the IMA Journal of Mathematical Control and Information. His current research interests include nonlinear control theory, hybrid systems, and control of partial differential equations, with applications including navigation and object tracking, fluid dynamics, and fusion control. He is an IMA Fellow, and an IEEE Fellow.
	\end{IEEEbiography}
	
		\begin{IEEEbiography}[{\includegraphics[width=1in,height=1.25in,clip,keepaspectratio]{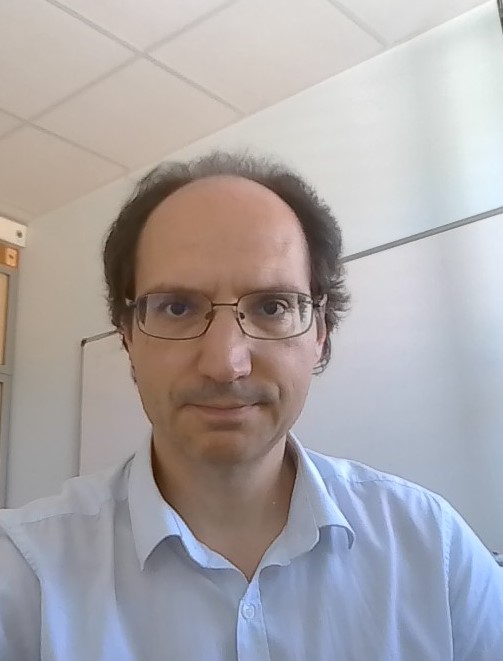}}]{Gildas Besançon}  received an Electrical Engineering diploma in 1993 and a PHD degree in Control in 1996, both in Grenoble, France. He became a full professor in 2010, and is currently affiliated to Grenoble Institute of Engineering and Management of Univ. Grenoble Alpes, France. He was also distinguished in 2010 by becoming a member of  Institut Universitaire de France. His research activities are held at the Control and Diagnosis
		Department of Gipsa-lab in the area of nonlinear and complex systems, with
		a special interest in observer issues, and applications in energy, hydraulics,
		micro/nano-sciences among other ones. He currently serves as an Associate Editor for IEEE Transactions of Automatic Control, as well as in various IEEE and IFAC Technical Committees.
	\end{IEEEbiography}
	
		\begin{IEEEbiography}[{\includegraphics[width=1in,height=1.25in,clip,keepaspectratio]{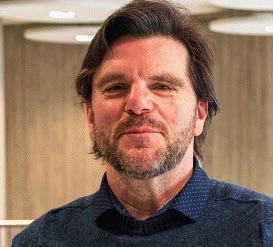}}]{Joseph J. Winkin} is full professor at the University of Namur, Belgium. He is also member of the Namur Institute for Complex Systems (naXys), University of Namur, Belgium. His main research interest is in the area of system and control theory, especially distributed parameter (infinite-dimensional) system theory and applications, linear-quadratic optimal control, spectral factorization techniques and dynamical analysis and control of (bio)chemical reactor models.
			He has also research interest in nonlinear system analysis and control and in application of control to pharmacokinetics. J.J. Winkin has contributed to numerous technical papers in these fields. He is member of the IEEE and IFAC Technical Committees on Distributed Parameter Systems. He is also member of IEEE, BMS (Belgian Mathematical Society) and EMS (European Mathematical Society), and serves as a reviewer for several journals and conferences in the system and control area. He is a past associate editor of the IEEE Transactions on Automatic Control.

				\end{IEEEbiography}
			
\end{document}